\documentclass[leqno,11pt]{amsart}

\usepackage[top=3.75cm, bottom=3cm, left=3cm, right=3cm]{geometry}
\usepackage[dvipsnames, table]{xcolor}
\usepackage{amsthm}
\theoremstyle{plain}
\usepackage{amssymb}
\usepackage{marvosym}
\usepackage{bm}
\usepackage{mathrsfs}
\usepackage{enumerate}  
\usepackage{mathtools}
\usepackage{tikz-cd}
\usepackage{graphicx}
\usepackage{float}
\usepackage[titletoc,title]{appendix}
\usepackage{marginnote}
\usepackage[disable]{todonotes}
\usepackage{etoolbox}
\usepackage[all]{xy}
\makeatletter
\patchcmd{\Ginclude@eps}{"#1"}{#1}{}{}
\makeatother
\usepackage[outdir=./]{epstopdf}
\usepackage[numbers]{natbib}
\usepackage[utf8]{inputenc}
\usepackage[english]{babel}
\usepackage{changepage}
\usepackage{makecell}
\usepackage{enumitem}

\definecolor{lightblue}{HTML}{1F88CD}
\definecolor{lightgrey}{HTML}{727272}
\definecolor{lightblue2}{HTML}{009EC1}
\definecolor{mypink}{HTML}{FD00B0}
\definecolor{lightred}{HTML}{ff4d4d}

\usepackage{hyperref}
\hypersetup{
	colorlinks=true,
    linkcolor={OliveGreen},
    citecolor={lightred},
	urlcolor={black}
}

\newtheorem*{theorem*}{Theorem}
\newtheorem{theorem}{Theorem}[section]
\newtheorem{corollary}[theorem]{Corollary}
\newtheorem{lemma}[theorem]{Lemma}
\newtheorem{conjecture}[theorem]{Conjecture}

\newtheorem{proposition}[theorem]{Proposition}
\theoremstyle{definition}
\newtheorem{example}[theorem]{Example}
\theoremstyle{definition}
\newtheorem{definition}[theorem]{Definition}
\theoremstyle{definition}
\newtheorem{remark}[theorem]{Remark}
\theoremstyle{definition}

\theoremstyle{definition}

\theoremstyle{definition}

\theoremstyle{definition}

\theoremstyle{definition}

\theoremstyle{definition}
\newtheorem{question!}[theorem]{Question!}
\theoremstyle{definition}

\makeatletter
\newcommand*\sbt{\mathpalette\sbt@{.75}}
\newcommand*\sbt@[2]{\mathbin{\vcenter{\hbox{\scalebox{#2}{$\m@th#1\bullet$}}}}}
\makeatother

\newcommand{\ra}{\rightarrow}

\newcommand{\sst}{\subset}

\newcommand{\bJ}{\bm{\mathrm{J}}}

\newcommand{\D}{\mathrm{D}}

\newcommand{\sH}{\mathrm{H}}

\newcommand{\ch}{\mathrm{ch}}

\newcommand{\HH}{\mathrm{HH}}

\newcommand{\dg}{\mathrm{dg}}

\newcommand{\HT}{\mathrm{HT}}

\DeclareMathOperator{\identity}{id}

\DeclareMathOperator{\Ext}{Ext}

\DeclareMathOperator{\Hom}{Hom}
\DeclareMathOperator{\RHom}{RHom}
\DeclareMathOperator{\Perf}{Perf}

\DeclareMathOperator{\Pic}{Pic}

\DeclareMathOperator{\Gr}{Gr}

\DeclareMathOperator{\HP}{\mathrm{HP}}

\newcommand{\cC}{\mathcal{C}}
\newcommand{\cA}{\mathcal{A}}
\newcommand{\cE}{\mathcal{E}}
\newcommand{\cU}{\mathcal{U}}

\newcommand{\cB}{\mathcal{B}}

\newcommand{\cQ}{\mathcal{Q}}
\newcommand{\Ku}{\mathcal{K}u}
\newcommand{\cP}{\mathcal{P}}
\newcommand{\cD}{\mathcal{D}}

\newcommand{\cN}{\mathcal{N}}
\newcommand{\cM}{\mathcal{M}}

\DeclareMathOperator{\oh}{\mathcal{O}}

\usetikzlibrary{decorations.pathmorphing}

\title[Infinitesimal categorical Torelli theorems for Fano threefolds]{Infinitesimal categorical Torelli theorems for Fano threefolds} 


\address{School of Mathematics, The University of Edinburgh, James Clerk Maxwell Building, Kings Buildings, Edinburgh, United Kingdom, EH9 3FD}
\email{a.jacovskis@sms.ed.ac.uk}

\address{Yau Mathematical Sciences Center, Tsinghua University, Beijing, China}
\email{lin-x18@mails.tsinghua.edu.cn}

\address{Institute for Advanced Study in Mathematics, Zhejiang University, Hangzhou, Zhejiang Province 310030, P. R. China}
\address{College of Mathematics, Sichuan University, Chengdu, Sichuan Province 610064, P. R. China}
\email{zhiyuliu@stu.scu.edu.cn, jasonlzy0617@gmail.com}
\urladdr{sites.google.com/view/zhiyuliu}

\address{Max Planck Institute for Mathematics, Vivatsgasse 7, 53111 Bonn, Germany}
\address{Institut de Mathématiqes de Toulouse, UMR 5219, Université de Toulouse, Université Paul Sabatier, 118 route de
Narbonne, 31062 Toulouse Cedex 9, France}
\email{shizhuozhang@mpim-bonn.mpg.de,shizhuo.zhang@math.univ-toulouse.fr}

\author{Augustinas Jacovskis, Xun Lin, Zhiyu Liu and Shizhuo Zhang}
\date{\today}

\begin{document}

\begin{abstract}
Let $X$ be a smooth Fano variety and $\Ku(X)$ its Kuznetsov component. A Torelli theorem for $\Ku(X)$ states that $\Ku(X)$ is uniquely determined by a certain polarized abelian variety associated to it. An infinitesimal Torelli theorem for $X$ states that the differential of the period map is injective. A categorical variant of the infinitesimal Torelli theorem for $X$ states that the morphism $\eta : H^1(X,T_X) \rightarrow \mathrm{HH}^2(\Ku(X))$ is injective. In the present article, we use the machinery of Hochschild (co)homology to relate the aforementioned three Torelli-type theorems for smooth Fano varieties via a commutative diagram. As an application, we prove infinitesimal categorical Torelli theorems for a class of prime Fano threefolds. We then prove, infinitesimally, a restatement of the Debarre--Iliev--Manivel conjecture regarding the general fiber of the period map for ordinary Gushel--Mukai threefolds. 
\end{abstract}


\maketitle

{
\hypersetup{linkcolor=black}
}

\section{Introduction}
Torelli problems are some of the oldest and the most classical problems in various aspects of algebraic geometry, including Hodge theory, birational geometry, moduli spaces of algebraic varieties, etc. The classical Torelli question asks whether an algebraic variety $X$ is uniquely determined by an abelian variety associated to it. Denote by $\mathcal{P}$ the period map $\mathcal{P}:\mathcal{M}\rightarrow\mathcal{D}$, where $\mathcal{M}$ is the moduli space of some class of algebraic varieties up to isomorphism, and $\mathcal{D}$ is the period domain. A \emph{Torelli theorem} holds for $X$ if $\mathcal{P}$ is injective. An \emph{infinitesimal Torelli theorem} holds for $X$ if the differential $d\mathcal{P}$ of the period map is injective. If $X$ is a smooth Fano threefold of Picard rank $1$ -- the focus of the second part of our paper -- then the period map $\mathcal{P}$ is given by $X\mapsto J(X)$, where $J(X)$ is the intermediate Jacobian of $X$. 

On the other hand, the seminal work \cite{bondal2001reconstruction} shows that the bounded derived category $\D^b(X)$ of a smooth projective Fano variety determines the isomorphism class of $X$. In other words, a \emph{derived Torelli theorem} holds for $X$. It is natural to ask for a class of Fano varieties whether they are also determined by \emph{less} information than the whole derived category. 

A natural candidate is a subcategory $\Ku(X)\subset \D^b(X)$ called the \emph{Kuznetsov component}, which is defined as the semiorthogonal complement of a natural exceptional collection of vector bundles on $X$. It is widely believed that the Kuznetsov component encodes essential geometric information of $X$, and it has been studied extensively by Kuznetsov and others (e.g. in \cite{kuznetsov2003derived, kuznetsov2009derived, kuznetsov2018derived}) for many Fano varieties. In particular, for a smooth cubic threefold $X$, its Kuznetsov component $\Ku(X)$ determines its isomorphism class (see \cite{bernardara2012categorical} and \cite{pertusi2020some}). This is known in the literature as a \emph{categorical Torelli theorem}. 

As in the case of the classical Torelli problem, one could imagine that the association $X\mapsto\Ku(X)$ is a ``categorical period map" $\mathcal{P}_{\mathrm{cat}}$ lifting the classical period map, defined on the moduli space $\mathcal{M}$ of smooth Fano threefolds up to isomorphism.  But since there is no good notion of a moduli space of triangulated categories (or dg categories) $\Ku(X)\subset\D^b(X)$, we cannot make sense of $\mathcal{P}_{\mathrm{cat}}$ mathematically. Nevertheless, the map $\eta \colon H^1(X,T_X)\rightarrow \HH^2(\mathcal{K}u(X))$ can be interpreted as its differential. We say that an \emph{infinitesimal categorical Torelli theorem} holds for $X$ if the map $\eta$ is injective. 

Recently, the intermediate Jacobian of a smooth Fano variety $X$ was reconstructed from its Kuznetsov component $\Ku(X)$ in \cite{perry2020integral}. This relates the Hodge-theoretic and categorical invariants of $X$ in one direction, in the sense that an equivalence of Kuznetsov components implies an isomorphism of intermediate Jacobians.



\subsection{Main Results}

\subsubsection{Infinitesimal Torelli vs. infinitesimal categorical Torelli}


In the present article, we relate infinitesimal Torelli theorems and infinitesimal categorical Torelli theorems for a class of Fano threefolds of Picard rank one, using the machinery of Hochschild (co)homology. We write $$\sH\Omega_{\bullet}(X)=\bigoplus_{p-q=\bullet}H^{p}(X,\Omega_{X}^{q})$$
 and
 $$\HT^{\bullet}(X)=\bigoplus_{p+q=\bullet}H^{p}(X,\Lambda^{q}T_{X}).$$
We prove the following theorem:

\begin{theorem}[Theorem \ref{commutative diagram theorem}] \label{intro infinitesimal diagram}
Let $X$ be a smooth projective variety. Assume there is a semiorthogonal decomposition $\D^b(X) = \langle \Ku(X), E_1, \dots, E_n \rangle$, where $\{ E_1, \dots, E_n \}$ is an exceptional collection. Then we have a commutative diagram
$$\xymatrix@C8pc@R2pc{\HH^{2}(\Ku(X))\ar[r]^{\gamma}&\Hom(\sH\Omega_{-1}(X),\sH\Omega_{1}(X))\\ \HT^{2}(X)\ar[u]^{\alpha'}\ar[ru]^{\tau}&\\
H^{1}(X,T_{X})\ar[ruu]^{\mathrm{d} \cP}\ar[u]^{\mathrm{inclusion}}&}$$
where $\tau$ is defined as a contraction of polyvector fields.
\end{theorem}

\begin{remark}
The map $\eta \colon H^{1}(X,T_{X})\longrightarrow \HH^{2}(\Ku(X))$ is defined as the composition of the vertical maps in the commutative diagram. 
\end{remark}

If we assume that the Hochschild homology of $X$ vanishes for degrees greater than one, the space of first-order deformations of $J(\Ku(X))$ is given by $\Hom(\HH_{-1}(\Ku(X)),\HH_{1}(\Ku(X)))$ by Proposition \ref{deformationtorus}. If we assume further that $J(\Ku(X))\cong J(X)$ as abelian varieties, then the diagram from the theorem above can be interpreted as taking tangent spaces of the diagram
\[
\begin{tikzcd}
\{ \Ku(X) \} / \simeq \arrow[r]         & \text{\{Abelian Varieties \}} / \cong \\
\{X\}/\cong \arrow[u] \arrow[ru] &         
\end{tikzcd}
\]
provided we have good knowledge of the moduli spaces in question.


\begin{corollary}\label{classicalandcategoricalinfinitesimal}
Infinitesimal classical Torelli for $X$ implies infinitesimal categorical Torelli for the Kuznetsov component $\Ku(X)$.
\end{corollary}

\begin{proof}
Suppose $\mathrm{d} \cP$ is injective. Then the fact that we have a composition $\mathrm{d} \cP = \gamma \circ \eta$ implies that $\eta$ is injective too.
\end{proof}

The main examples of infinitesimal categorical Torelli theorems that we study are those for Picard rank one, index one and two Fano threefolds. Recall that Fano threefolds satisfy the assumption $\HH_{2i+1}(X) = 0$ for $i\geq 1$ and $J(\Ku(X))\cong J(X)$ as abelian varieties (see 
Lemma~\ref{assumption}). We summarise our results in the following theorem.

\begin{theorem}[Theorems \ref{Y1234}, \ref{Xo10}, \ref{X12141618}, \ref{X248}]
Let $X_{2g-2}$ be a Fano threefold of index one and degree $2g-2$, where $g$ is its genus. Let $Y_d$ be a Fano threefold of index two and degree $d$. 
\begin{itemize}
    \item For $Y_d$ where $1 \leq d \leq 4$, the infinitesimal categorical Torelli theorem holds.
    \item For $X_{2g-2}$ where  $g=2,4,5,7$, the infinitesimal categorical Torelli theorem holds and the same result holds for $X_4$ if it is not hyperelliptic. 
\end{itemize}
\end{theorem}

\subsubsection{The Debarre--Iliev--Manivel Conjecture}

In \cite{debarre2012period}, the authors conjecture that the general fiber of the classical period map from the moduli space of ordinary Gushel--Mukai (GM) threefolds to the moduli space of $10$-dimensional principally polarised abelian varieties is the disjoint union of $\cC_m(X)$ (the minimal surface of the Hilbert scheme of conics on $X$) and $M_G^X(2,1,5)$ (a Gieseker moduli space of semistable sheaves of rank two, $c_1=H, c_2=5L$ on $X$), both quotiented by involutions. We call this the \emph{Debarre--Iliev--Manivel Conjecture}.

Within the moduli space of smooth GM threefolds, we define the fiber of the ``categorical period map" at $\Ku(X)$ as the set of isomorphism classes of all ordinary GM threefolds $X'$ whose Kuznetsov components satisfy $\Ku(X')\simeq\Ku(X)$. In our recent work \cite{JLZ2021torelli}, we prove the categorical analogue of the Debarre--Iliev--Manivel conjecture:

\begin{theorem}[{\cite[Theorem 1.7]{JLZ2021torelli}}] \label{intro categorical period map}
A general fiber of the categorical period map at the Kuznetsov component $\Ku(X)$ of an ordinary GM threefold $X$ is the union of $\mathcal{C}_m(X)/\iota$ and $M_G^X(2,1,5)/\iota'$, where $\iota, \iota'$ are geometrically meaningful involutions.
\end{theorem}

As an application, the Debarre--Iliev--Manivel conjecture can be restated in an equivalent form as follows:

\begin{conjecture}
\label{DIM_equi_conjecture}
Let $X$ be a general ordinary GM threefold. Then the intermediate Jacobian $J(X)$ determines the Kuznetsov component $\Ku(X)$.
\end{conjecture}

Although we are not able to prove the conjecture, we are able to show that an infinitesimal version of it holds:

\begin{theorem}[{Theorem \ref{Xo10}}]
\label{DIM_conjecture_infinitesimal}
Let $X$ be an ordinary GM threefold. Then the map \[ \gamma \colon \HH^2(\Ku(X)) \longrightarrow \Hom(\sH\Omega_{-1}(X), \sH\Omega_1(X)) \] is injective. 
\end{theorem}





\subsection{Organization of the paper}
In Section \ref{sod section}, we collect basic facts about semiorthogonal decompositions. In Section \ref{hochschild cohomology section}, we recall the definition of Hochschild (co)homology for admissible subcategories of bounded derived categories $\D^b(X)$ of smooth projective varieties $X$. We then prove Theorem~\ref{intro infinitesimal diagram}. In Section \ref{infinitesimal section}, we apply the techniques developed in Section \ref{hochschild cohomology section} to prime Fano threefolds of index one and two. In particular, we show the infinitesimal version (Theorem \ref{DIM_conjecture_infinitesimal}) of Conjecture~\ref{DIM_equi_conjecture} for ordinary GM threefolds.

\subsection*{Acknowledgements}
Firstly, it is our pleasure to thank Arend Bayer for very useful discussions on the topics of this project. We would like to thank Sasha Kuznetsov for answering many of our questions on Gushel--Mukai threefolds. We thank  Enrico Fatighenti and Luigi Martinelli for helpful conversations on several related topics. We thank the anonymous referees for their careful reading of our manuscript, and their many insightful comments and suggestions. The first and last authors were supported by ERC Consolidator Grant WallCrossAG, no. 819864. The last author is also supported by ANR project FanoHK, grant ANR-20-CE40-0023. The third author is partially supported by NSFC (nos. 11890660 and 11890663). Part of this work was finished while the last author was visiting the Max-Planck Institute For Mathematics. He is grateful for their excellent hospitality and support.

\section{Semiorthogonal decompositions} \label{sod section}

In this section, we collect some useful facts/results about semiorthogonal decompositions. Background on triangulated categories and derived categories of coherent sheaves can be found in \cite{huybrechts2006fourier}, for example. From now on, let $\D^b(X)$ denote the bounded derived category of coherent sheaves on a smooth projective variety $X$. We use $k$ to denote the field of complex numbers $\mathbb{C}$.


\subsection{Exceptional collections and semiorthogonal decompositions}

\begin{definition} Let $\cD$ be a triangulated category. We say that $E \in \cD$ is an \emph{exceptional object} if $\RHom(E, E) = k$. Now let $(E_1, \dots, E_m)$ be a collection of exceptional objects in $\cD$. We say it is an \emph{exceptional collection} if $\RHom(E_i, E_j) = 0$ for $i > j$.
\end{definition}

\begin{definition}
Let $\cD$ be a triangulated category and $\cC$ a triangulated subcategory. We define the \emph{right orthogonal complement} of $\cC$ in $\cD$ as the full triangulated subcategory
\[ \cC^\bot = \{ X \in \cD \mid \Hom(Y, X) =0  \text{ for all } Y \in \cC \}.  \]
The \emph{left orthogonal complement} is defined similarly, as 
\[ {}^\bot \cC = \{ X \in \cD \mid \Hom(X, Y) =0  \text{ for all } Y \in \cC \}.  \]
\end{definition}

\begin{definition}
Let $\cD$ be a triangulated category. We say a triangulated subcategory $\cC \sst \cD$ is \emph{admissible}, if the inclusion functor $i : \cC \hookrightarrow \cD$ has left adjoint $i^*$ and right adjoint $i^!$.
\end{definition}

\begin{definition}
Let $\cD$ be a triangulated category, and $( \cC_1, \dots, \cC_m  )$ be a collection of full admissible subcategories of $\cD$. We say that $\cD = \langle \cC_1, \dots, \cC_m \rangle$ is a \emph{semiorthogonal decomposition} of $\cD$ if $\cC_j \sst \cC_i^\bot $ for all $i > j$, and the subcategories $(\cC_1, \dots, \cC_m )$ generate $\cD$, i.e. the smallest strictly full triangulated subcategory of $\cD$ containing $\cC_1,\dots,\cC_m$ is equal to $\cD$. 
\end{definition}


\begin{definition}
The \emph{Serre functor} $S_{\cD}$ of a triangulated category $\cD$ is the autoequivalence of $\cD$ such that there is a functorial isomorphism of vector spaces
\[ \Hom_{\cD}(A, B) \cong \Hom_{\cD}(B, S_{\cD}(A))^\vee \]
for any $A, B \in \cD$.
\end{definition}

\begin{proposition} \label{serre s.o.d. proposition}
Assume a triangulated category $\cD$ admits a Serre functor $S_{\cD}$ and let $\cD = \langle \cD_1 , \cD_2 \rangle$ be a semiorthogonal decomposition. Then $\cD \simeq \langle S_{\cD}(\cD_2), \cD_1 \rangle \simeq \langle \cD_2, S^{-1}_{\cD}(\cD_1) \rangle$ are also semiorthogonal decompositions.
\end{proposition}

\begin{example}
Let $X$ be a smooth projective variety and $\cD = \D^b(X)$. Then $S_X:=S_{\cD}(-) = (- \otimes \mathcal{O}(K_X))[\dim X]$.
\end{example}

\section{Hochschild (co)homology and infinitesimal Torelli theorems} \label{hochschild cohomology section}

\subsection{Definitions}\label{sectionHochschild}

In this subsection, we recall some basics on Hochschild (co)homology of admissible subcategories of $\D^{b}(X)$, where $X$ is a smooth projective variety. We refer the reader to \cite{kuznetsov2009hochschild} for more details. For background on Hochschild (co)homology of dg-categories, we refer the reader to the paper \cite{KELLER1998223} and survey \cite{Keller2006OnDG}.

\begin{definition}\label{Hochschilddg}
 Let $\mathcal{C}$ be a dg-category. The Hochschild cohomology of $\mathcal{C}$ is 
 $$\HH^{\bullet}(\mathcal{C}) \coloneqq \Hom_{\D(\mathcal{C}\otimes \mathcal{C}^{\mathrm{opp}})}(\mathcal{C},\mathcal{C}[\bullet]).$$
 The Hochschild homology of $\mathcal{C}$ is 
 $$\HH_{\bullet}(\mathcal{C}):=H^{\bullet}(\mathcal{C}\otimes^{\mathbb{L}}_{\mathcal{C}\otimes\mathcal{C}^{\mathrm{opp}}}\mathcal{C}).$$

 We write $\bullet$ to represent integers. Note that the degree $i$ Hochschild homology of $\mathcal{C}$ is usually defined as $H^{-i}(\mathcal{C}\otimes^{\mathbb{L}}_{\mathcal{C}\otimes\mathcal{C}^{\mathrm{opp}}}\mathcal{C})$ (=$H_{i}(\mathcal{C}\otimes^{\mathbb{L}}_{\mathcal{C}\otimes\mathcal{C}^{\mathrm{opp}}}\mathcal{C})$) in the standard literature/context. The sign convention that we use is to make our paper consistent with Hochschild (co)homology of admissible subcategories in the paper \cite{kuznetsov2009hochschild}. The object $\bigoplus_{i} \HH_{i}(\mathcal{C})$ is a module over the algebra $\bigoplus_{j}\HH^{j}(\mathcal{C})$ defined by compositions. Using our sign convention, the degree $j$ cohomology acts on degree $i$ homology to land in degree $i+j$ homology. Using the standard sign convention, the degree $j$ cohomology acts on degree $i$ homology to land in degree $i-j$ homology.
\end{definition}

\begin{definition}[{\cite{kuznetsov2009hochschild}}]
Let $X$ be a smooth projective variety, and $\mathcal{A}$ be an admissible subcategory of $\D^{b}(X)$. Consider any semiorthogonal decomposition of $\D^{b}(X)$ that contains $\mathcal{A}$ as a component. Let $P$ be the kernel of the projection to $\mathcal{A}$. The \emph{Hochschild cohomology} of $\mathcal{A}$ is defined as 
$$\HH^{\bullet}(\mathcal{A}) \coloneqq \Hom_{\D^{b}(X\times X)}(P,P[\bullet]).$$
The \emph{Hochschild homology} of $\mathcal{A}$ is defined as 
$$\HH_{\bullet}(\mathcal{A}) \coloneqq \Hom_{\D^{b}(X\times X)}(P,P\circ S_{X}[\bullet]).$$
\end{definition}

\begin{lemma}[{\cite[Theorem 4.5, Proposition 4.6]{kuznetsov2009hochschild}}] \label{dgalgebra} 
Let $E$ be a strong compact generator of $\D^{b}(X)$, and let $E_{\mathcal{A}}$ be the projection to $\mathcal{A}$. Then $E_{\mathcal{A}}$ is a strong compact generator of $\mathcal{A}$. Define $C =  \RHom(E_{\mathcal{A}},E_{\mathcal{A}})$. Then there are isomorphisms
$$\HH^{\bullet}(C)\cong \HH^{\bullet}(\mathcal{A}) \  \text{ and } \ \HH_{\bullet}(C)\cong \HH_{\bullet}(\mathcal{A}).$$
\end{lemma}

\begin{remark} \leavevmode 
Let $F$ be any strong compact generator of $\mathcal{A}$. Define $A=\RHom(F,F)$. 
\begin{enumerate}
    \item Let $\mathrm{Perf}_{\dg}(X)$ be a dg-enhancement of $\mathrm{Perf}(X)$ whose objects are K-injective perfect complexes, and let $\mathcal{A}_{\mathrm{dg}}$ be a dg-subcategory of $\mathrm{Perf}_{\mathrm{dg}}(X)$ whose objects are in $\mathcal{A}$. Then  we have the isomorphisms
    $$\HH_{\bullet}(A)\cong\HH_{\bullet}(\mathcal{A}_{\dg}) \ \text{ and } \ \HH^{\bullet}(A)\cong\HH^{\bullet}(\mathcal{A}_{\dg})$$
    because the morphism of dg-categories $A \ra \mathcal{A}_{\mathrm{dg}}$ is a derived Morita equivalence. Here $A$ is a dg-category with one object and the endomorphism of this unique object is the dg-algebra $A$; the morphism sends the unique object to an K-injective resolution of $F$.
    \item Since $\Perf(X)$ is a triangulated category with a unique enhancement \cite{Lunts_2010}, the Hochschild (co)homology of admissible subcategories of $\Perf(X)$ defined by Kuznetsov coincides with that of the subcategory of the dg-enhancement that naturally comes from the dg-enhancement of $\Perf(X)$.
\end{enumerate}
\end{remark}



Let $A$ be a $k$-algebra. Note that the Hochschild homology $\HH_{\bullet}(A)$ is a graded $\HH^{\bullet}(A)$-module. The module structure is easily described by the definition of Hochschild (co)homology via the $\Ext$ and $\mathrm{Tor}$ functors. Consider a variety $X$. By the Hochschild--Kostant--Rosenberg (HKR) isomorphism, the degree two Hochschild cohomology has a summand $H^{1}(X,T_{X})$ which is the first-order deformations of $X$ as a variety. The action of $H^{1}(X,T_{X})$ on the Hochschild homology via the module structure can be interpreted as deformations of a certain invariant with respect to the deformations of $X$. Here, we have the invariant $\HH_{\bullet}(X)\cong\bigoplus_{p-q=\bullet}H^{p}(X,\Omega^{q}_{X})$ that is closely related to the intermediate Jacobian of $X$. 


In the case of admissible subcategories of derived categories, we can describe the module structure by kernels. For more background on the calculus of kernels, see \cite{kuznetsov2009hochschild}.

\begin{definition}\label{moduleaction}
Let $\mathcal{A}$ be an admissible subcategory of $\D^{b}(X)$, and let $P$ be the kernel of the left projection to $\mathcal{A}$. Denote by $S_X = \Delta_{*} \omega_X [ \mathrm{dim} X ]$ the kernel of the Serre functor of $\D^b(X)$.
Take $\alpha\in \HH^{a}(\mathcal{A})$ and $\beta\in \HH_{b}(\mathcal{A})$. The action of $\alpha$ on $\beta$ is the composition 
\[ P \xrightarrow{\, \, \, \beta \, \, \,} P \circ S_X[b] \xrightarrow{\, \alpha \otimes \identity \, } P \circ S_X[a+b] . \]
We also have that the object $\bigoplus_{i}\Hom_{\D^{b}(X\times X)}(P,P\circ S_{X}[i])$ is a module over the graded algebra $\bigoplus_{j}\Hom_{\D^{b}(X\times X)}(P,P[j])$ by this action.
\end{definition}

\begin{proposition}
Let $\mathcal{A}$ be an admissible subcategory of $\Perf(X)$. Let $E$ be a strong compact generator of $\mathcal{A}$, and $A:= \RHom(E,E)$. The isomorphisms $\HH^{\bullet}(A)\cong \HH^{\bullet}(\mathcal{\mathcal{A}})$ and $\HH_{\bullet}(A)\cong \HH_{\bullet}(\mathcal{A})$ from Lemma \ref{dgalgebra} preserve both sides of the obvious module structure and algebra structure of Hochschild cohomology.\todo[color=lightgrey!25]{Changed "preserve two sides obvious module structure and algebra structure of Hochschild cohomology." to what it is now.}
\end{proposition}

\begin{proof}
We follow \cite[Theorem 4.5, Proposition 4.6]{kuznetsov2009hochschild}. Consider the semiorthogonal decomposition $\Perf(X)=\langle \mathcal{A}, {}^\bot \mathcal{A}\rangle$.
By \cite[Theorem 3.1]{ARMENTA2019236}, we may assume 
the compact generator $E$ is the one in \cite[Theorem 4.5]{kuznetsov2009hochschild}. Let $F$ be the compact generator of $\mathcal{B}^{\vee}$ introduced in \cite[Theorem 4.5]{kuznetsov2009hochschild}. See \emph{loc. cit.} for the definition of $\cB^\vee$.
There is a fully faithful functor $$\mu: \Perf(A\otimes A^{\mathrm{opp}})\rightarrow \D^{b}(X\times X)$$ such that $\mu(A)=P$. Explicitly $P=A\otimes^{\mathbb{L}}_{A\otimes A^{\mathrm{opp}}}E\boxtimes F$. Thus $\Hom_{\Perf(A\otimes A^{\mathrm{opp}})}(A,A[\bullet])\cong \Hom_{\D^{b}(X\times X)}(P,P[\bullet])$, which is compatible with the algebra structure since both of the algebra structures are defined by compositions. It remains to check the compatibility of the module structure.

\par

There are isomorphisms \cite[Theorem 4.5]{kuznetsov2009hochschild}
$$\mathrm{R}\Gamma(X\times X,(E\boxtimes F)\otimes^{\mathbb{L}} (F\boxtimes E))\cong \mathrm{R}\Gamma(X,E\otimes^{\mathbb{L}} F)\otimes \mathrm{R}\Gamma(X,E\otimes^{\mathbb{L}} F)\cong A\otimes A^{\mathrm{opp}}.$$
Therefore we have a morphism 
$$\xymatrix@C=0.5cm@R=0.5cm{A\otimes^{\mathbb{L}}_{A\otimes A^{\mathrm{opp}}}(A\otimes A^{\mathrm{opp}})\otimes^{\mathbb{L}}_{A\otimes A^{\mathrm{opp}}} A\ar[r]&
\mathrm{R}\Gamma(X\times X, A\otimes^{\mathbb{L}}_{A\otimes A^{\mathrm{opp}}}(E\boxtimes F)\otimes^{\mathbb{L}} (F\boxtimes E)\otimes^{\mathbb{L}}_{A^{\mathrm{opp}}\otimes A}A)\\
A\otimes^{\mathbb{L}}_{A\otimes A^{\mathrm{opp}}}A\ar@{=}[u]&
},$$
which is a quasi-isomorphism by dévissage. Namely, by using a semi-free resolution, it suffices to check the case $A\otimes A^{\mathrm{opp}}$. The quasi-isomorphism is compatible with the module structure. Namely, for an element $\alpha\in \Hom_{\Perf(A\otimes A^{\mathrm{opp}})}(A,A[j])$, we have a commutative diagram
$$\xymatrix{A\otimes^{\mathbb{L}}_{A\otimes A^{opp}}A\ar[r]\ar[d]^{\alpha}&\mathrm{R}\Gamma(X\times X, A\otimes^{\mathbb{L}}_{A\otimes A^{\mathrm{opp}}}(E\boxtimes F)\otimes^{\mathbb{L}} (F\boxtimes E)\otimes^{\mathbb{L}}_{A^{\mathrm{opp}}\otimes A}A)\ar[d]^{\alpha}\\
A[j]\otimes^{\mathbb{L}}_{A\otimes A^{opp}}A\ar[r]&\mathrm{R}\Gamma(X\times X, A[j]\otimes^{\mathbb{L}}_{A\otimes A^{\mathrm{opp}}}(E\boxtimes F)\otimes^{\mathbb{L}} (F\boxtimes E)\otimes^{\mathbb{L}}_{A^{\mathrm{opp}}\otimes A}A)}$$
Then for any integer $i$ we have the commutative diagram
$$\xymatrix{H^{i}(A\otimes^{\mathbb{L}}_{A\otimes A^{opp}}A)\ar[r]\ar[d]^{\alpha}&H^{i}(X\times X, A\otimes^{\mathbb{L}}_{A\otimes A^{\mathrm{opp}}}(E\boxtimes F)\otimes^{\mathbb{L}} (F\boxtimes E)\otimes^{\mathbb{L}}_{A^{\mathrm{opp}}\otimes A}A)\ar[d]^{\alpha}\\
H^{i+j}(A\otimes^{\mathbb{L}}_{A\otimes A^{opp}}A)\ar[r]&H^{i+j}(X\times X, A\otimes^{\mathbb{L}}_{A\otimes A^{\mathrm{opp}}}(E\boxtimes F)\otimes^{\mathbb{L}} (F\boxtimes E)\otimes^{\mathbb{L}}_{A^{\mathrm{opp}}\otimes A}A)}$$
Finally, we must check that the isomorphism $H^{\bullet}(X\times X, P\otimes^{\mathbb{L}} P^{T}) \cong \Hom_{\D^{b}(X\times X)}(P,P\circ S_{X}[\bullet])$\footnote{$T$ is a transposition of the product $X\times X$. Namely, it maps the point $(x_{1},x_{2})$ to $(x_{2},x_{1})$.} in \cite[Proposition 4.6]{kuznetsov2009hochschild} is compatible with the module structure. Firstly, we have a functorial isomorphism with respect to the factor $P$:
\begin{align*}
    \Hom_{\D^{b}(X\times X)}(\mathcal{O}_{X\times X}, P\otimes^{\mathbb{L}} \Delta_{\ast}\mathcal{O}_{X}[\bullet])\cong& \Hom_{\D^{b}(X\times X)}((\Delta_{\ast}\mathcal{O}_{X})^{\vee}, P[\bullet])\\
    \cong& \Hom_{\D^{b}(X\times X)}((\Delta_{\ast}\mathcal{O}_{X})^{\vee}\circ S_{X},P\circ S_{X}[\bullet]).
\end{align*}
By Grothendieck Duality, $(\Delta_{\ast}\mathcal{O}_{X})^{\vee}\circ S_{X}\cong \Delta_{\ast}\mathcal{O}_{X}$, therefore
$$\Hom_{\D^{b}(X\times X)}((\Delta_{\ast}\mathcal{O}_{X})^{\vee}\circ S_{X},P\circ S_{X}[\bullet])\cong \Hom_{\D^{b}(X\times X)}(\Delta_{\ast}\mathcal{O}_{X},P\circ S_{X}[\bullet]).$$
Thus we have the following functorial isomorphism with respect to $P$:
$$\Hom_{\D^{b}(X\times X)}(\mathcal{O}_{X\times X}, P\otimes^{\mathbb{L}} \Delta_{\ast}\mathcal{O}_{X}[\bullet])\cong \Hom_{\D^{b}(X\times X)}(\Delta_{\ast}\mathcal{O}_{X},P\circ S_{X}[\bullet]).$$
Consider the triangle $P' \to \Delta_* \oh_X \to P$ with respect to the semiorthogonal decomposition $\langle \mathcal{A}, {}^\bot \mathcal{A}\rangle$. Note that $(\Delta_{\ast}\mathcal{O}_{X})^{T}=\Delta_{\ast}\mathcal{O}_{X}$, so we have a triangle $P'^{T} \to \Delta_* \oh_X \to P^{T}$. Tensoring by $P$, we get the triangle 
$P\otimes^{\mathbb{L}}P'^{T} \to P\otimes^{\mathbb{L}}\Delta_* \oh_X \to P\otimes^{\mathbb{L}}P^{T}$. Taking any element $\alpha\in \Hom_{\D^{b}(X\times X)}(P,P[j])$, we have a commutative diagram
$$\xymatrix{P\otimes^{\mathbb{L}}P'^{T}\ar[r]\ar[d]^{\alpha}&P\otimes^{\mathbb{L}}\Delta_{*}\oh_{X}\ar[r]\ar[d]^{\alpha}&P\otimes^{\mathbb{L}}P^{T}\ar[d]^{\alpha}\\
 P[j]\otimes^{\mathbb{L}}P'^{T}\ar[r]&P[j]\otimes^{\mathbb{L}}\Delta_{*}\oh_{X}\ar[r]&P[j]\otimes^{\mathbb{L}}P^{T}}$$
Applying $\Hom_{\D^{b}(X\times X)}(\mathcal{O}_{X\times X},-)$ to the commutative diagram above, we have the following commutative diagram for any integer $i$:
$$\xymatrix{\Hom_{\D^{b}(X\times X)}(\oh_{X\times X},P\otimes^{\mathbb{L}} \Delta_* \oh_X [i])\ar[r]^{\cong}\ar[d]^{\alpha}& \Hom_{\D^{b}(X\times X)}(\oh_{X\times X},P\otimes^{\mathbb{L}} P^{T}[i])\ar[d]^{\alpha}\\
  \Hom_{\D^{b}(X\times X)}(\oh_{X\times X},P\otimes^{\mathbb{L}} \Delta_* \oh_X [i+j])\ar[r]^{\cong}& \Hom_{\D^{b}(X\times X)}(\oh_{X\times X},P\otimes^{\mathbb{L}} P^{T}[i+j]) }$$
  The middle maps are isomorphisms since $\Hom_{\D^{b}(X\times X)}(\oh_{X\times X},P\otimes^{\mathbb{L}}P'^{T}[\bullet])=0$ (see \cite[Corollary 3.10]{kuznetsov2009hochschild}). Similarly, applying $\Hom_{\D^{b}(X\times X)}(-,P\circ S_{X})$ to the triangle $P' \to \Delta_* \oh_X \to P$, we get the commutative diagram 
$$\xymatrix@C=2cm{\Hom_{\D^{b}(X\times X)}(P,P\circ S_{X}[i])\ar[r]^{\cong}\ar[d]^{\alpha}&\Hom_{\D^{b}(X\times X)}(\Delta_{*}\oh_{X}, P\circ S_{X}[i])\ar[d]^{\alpha}\\
 \Hom_{\D^{b}(X\times X)}(P,P\circ S_{X}[i+j])\ar[r]^{\cong}&\Hom_{\D^{b}(X\times X)}(\Delta_{*}\oh_{X}, P\circ S_{X}[i+j]) }$$
since $\Hom_{\D^{b}(X\times X)}(P',P\circ S_{X}[\bullet])=0$ (see \cite[Cor 3.10]{kuznetsov2009hochschild}). Here $\alpha$ is the morphism induced by the morphism $\xymatrix{P\circ S_{X}\ar[r]^{\alpha}&P[j]\circ S_{X}}$.  

Combining the isomorphisms above, we get an isomorphism $$H^{\bullet}(X\times X, P\otimes^{\mathbb{L}} P^{T})\cong \Hom_{\D^{b}(X\times X)}(P,P\circ S_{X}[\bullet])$$ 
which is the isomorphism in \cite[Proposition 4.6]{kuznetsov2009hochschild}. It is compatible with the module structure over $\Hom_{\D^{b}(X\times X)}(P,P[\bullet])$ since each step preserves the module structure. Namely, for any element $\beta\in \Hom(P,P[j])$ and any integer $i$, there is a commutative diagram
$$\xymatrix{H^{i}(X\times X, P\otimes^{\mathbb{L}}P^{T})\ar[r]^{\cong}\ar[d]&\Hom(P,P\circ S_{X}[i])\ar[d]\\
H^{i+j}(X\times X, P\otimes^{\mathbb{L}}P^{T})\ar[r]^{\cong}&\Hom(P,P\circ S_{X}[i+j])}$$
The left morphism is induced by tensoring $\beta\in \Hom(P,P[j])$ with $P^T$, i.e. we have $$\xymatrix{P\otimes^{\mathbb{L}}P^{T}\ar[r]^{\beta}&P[j]\otimes^{\mathbb{L}}P^{T}}.$$
The right morphism is the map defined in Definition \ref{moduleaction}.
\end{proof}







\begin{theorem}\label{derivedinvariant}
Let $\mathcal{A}$ be an admissible subcategory of $\D^{b}(X)$, and $\mathcal{B}$ be an admissible subcategory of $\D^{b}(Y)$. Suppose the Fourier--Mukai functor $\Phi_{\mathcal{E}}: \D^{b}(X)\rightarrow \D^{b}(Y)$ induces an equivalence between the subcategories $\mathcal{A}$ and $\mathcal{B}$. Then we have isomorphisms of Hochschild cohomology 
$$\HH^{\bullet}(\mathcal{A})\cong \HH^{\bullet}(\mathcal{B})$$
and Hochschild homology
$$\HH_{\bullet}(\mathcal{A})\cong\HH_{\bullet}(\mathcal{B})$$ 
which preserve both sides of the module structure and algebra structure.
\end{theorem}

\begin{proof}
We have $\mathcal{A}_{\dg}\simeq\mathcal{B}_{\dg}$ in the homotopy category whose weak equivalences are Morita equivalences \cite[Proposition 9.4]{bernardara2014semiorthogonal}. Hence, there is an isomorphism $\Perf_{\dg}(\mathcal{A}_{\dg})\simeq \Perf_{\dg}(\mathcal{B}_{\dg})$ in $\mathsf{Hqe}$ \cite{10.1155/IMRN.2005.3309}. That is, $\mathcal{A}_{\dg}$ and $\mathcal{B}_{\dg}$ are connected by a chain of Morita equivalences
  \[
\begin{tikzcd}[column sep=0.8em]
                   & \Perf_\dg(\cA_\dg) \arrow[rd] &     & \cdots \arrow[ld] \arrow[rd] &                & \Perf_\dg(\cB_\dg) &                    \\
\cA_\dg \arrow[ru] &                              & C_0 & \quad \quad \cdots \quad   \quad                   & C_n \arrow[ru] &                   & \cB_\dg \arrow[lu]
\end{tikzcd}
 \]

By \cite[Theorem 3.1]{ARMENTA2019236}, if two dg-categories are derived equivalent induced by a bi-module (Morita equivalence), then the equivalence induces an isomorphism of Hochschild (co)homology and preserves the module structure.
\end{proof}

Let $X$ be a smooth algebraic variety. Classically we have the HKR isomorphisms \cite[Theorem 8.3]{kuznetsov2009hochschild} given by 
$$\mathrm{Hom}_{\D^{b}(X\times X)}(\mathcal{O}_{\Delta},\mathcal{O}_{\Delta}[\bullet])\cong \bigoplus_{p+q=\bullet} H^{p}(X,\wedge^{q}T_{X})$$
and
$$\mathrm{Hom}_{\D^{b}(X\times X)}(\mathcal{O}_{X\times X},\mathcal{O}_{\Delta}\otimes^{\mathbb{L}}\mathcal{O}_{\Delta}[\bullet])\cong \bigoplus_{p-q=\bullet}H^{p}(X,\Omega_{X}^{q}) . $$

However, the HKR isomorphisms may not preserve the obvious algebra structures and module structures. Let $\mathrm{IK}$ be the twist of the HKR isomorphisms with the square root of the Todd class. It was originally conjectured in \cite[Conjecture 5.2]{cualduararu2005mukai} and proved in \cite[Theorem 1.4]{calaque2012cualduararu} that $\mathrm{IK}$ is compatible with the module structures on differential forms over polyvector fields and on Hochschild homology over Hochschild cohomology.

\subsection{Deformations and infinitesimal (categorical) Torelli theorems}

In this subsection we will prove Theorem \ref{commutative diagram theorem} which relates the first-order deformations of a variety $X$, its intermediate Jacobian $J(X)$, and its Kuznetsov component $\Ku(X)$ via a commutative diagram.

We first recall the construction of the intermediate Jacobian of a triangulated subcategory of $\D^b(X)$.
Let $\mathcal{A}$ be an admissible subcategory of $\D^{b}(X)$. We write $\HP_{\bullet}(\mathcal{A})$ as periodic cyclic homology of $\mathcal{A}$. Denoted by $\mathrm{K}^{\mathrm{top}}_{\bullet}(\mathcal{A})$ is the topological $\mathrm{K}$-group of $\mathcal{A}$ defined by Blanc in \cite{blanc_2016}. Periodic cyclic homology and topological $\mathrm{K}$-theory are both $2$-periodic.
\begin{definition}[{\cite[Definition 5.24]{perry2020integral}}]
Let $\mathcal{A}$ be an admissible subcategory of $\D^{b}(X)$ and consider the diagram
  \[
\begin{tikzcd}
\mathrm{K}_1^\mathrm{top}(\cA) \arrow[r, "\ch_1^{\mathrm{top}}"] \arrow[rrd, "\pi'"'] & \HP_1(\cA) \arrow[r, "\cong"] & \bigoplus_n \HH_{2n-1}(\cA) \arrow[d, "\pi"]   \\
                                                                           &                               & \HH_1(\cA) \oplus \HH_3(\cA) \oplus \cdots
\end{tikzcd}
 \]
where $\pi$ is the natural projection, and $\pi'$ is the composition.
Define the \emph{intermediate Jacobian of $\cA$} as
 $$J(\mathcal{A})= (\HH_{1}(\mathcal{A})\oplus \HH_{3}(\mathcal{A})\oplus\cdots)/\Gamma,$$
 where $\Gamma$ is the image of $\pi'$. Note that $\Gamma$ is a lattice.
\end{definition}

\begin{lemma}\label{assumption}
 Let $X$ be a Fano threefold. Assume we have a semi-orthogonal decomposition $D^{b}(X)=\langle\Ku(X), E_{1}, E_{2}, \cdots, E_{n} \rangle$, where 
 $\langle E_{1}, E_{2}, \cdots, E_{n}\rangle$ is an exceptional collection. Then, $J(\Ku(X))\cong J(X)$
\end{lemma}

\begin{proof}
  The proof is essentially due to Alex Perry \cite{perry2020integral}.
  By \cite[Lemma 5.2]{perry2020integral}, there are Euler pairings
for $K^{top}_{-3}(\Ku(X))$ and $K^{top}_{-3}(D^{b}(X))$. Since $K^{top}_{-3}(\langle E\rangle)=0$ for any exceptional object $E$, we have 
$K^{top}_{-3}(\Ku(X))_{\mathrm{tf}}\cong K^{top}_{-3}(D^{b}(X))_{\mathrm{tf}}$ which preserves Euler pairings. Under assumptions of \cite[Proposition 5.23]{perry2020integral} (it is true for Fano threefolds), there is a Hodge isometry $K^{top}_{-3}(D^{b}(X))_{\mathrm{tf}}\cong H^{3}(X,\mathbb{Z})_{\mathrm{tf}}$, see the proof of \cite[Proposition 5.23]{perry2020integral}. Since the paring of $H^{3}(X,\mathbb{Z})_{\mathrm{tf}}$ is anti-symmetric, both of the Euler pairings of $K^{top}_{-3}(\Ku(X))_{\mathrm{tf}}$ and $K^{top}_{-3}(D^{b}(X))_{\mathrm{tf}}$ are anti-symmetric. Thus we have
a Hodge isometry $K^{top}_{-3}(\Ku(X))_{\mathrm{tf}}\cong H^{3}(X,\mathbb{Z})_{\mathrm{tf}}$, which implies an isomorphism of abelian varieties
$J(K^{top}_{-3}(\Ku(X)))\cong J(X)$. Since the topological $K$ group is $2$-periodic, we have $J(K^{top}_{-3}(\Ku(X)))\cong J(K^{top}_{1}(\Ku(X)):=J(\Ku(X))$, see \cite[Remark 5.25 ]{perry2020integral}.

\end{proof}
 




In the theorems below, we write $$\sH\Omega_{\bullet}(X)=\bigoplus_{p-q=\bullet}H^{p}(X,\Omega_{X}^{q})$$
 and
 $$\HT^{\bullet}(X)=\bigoplus_{p+q=\bullet}H^{p}(X,\Lambda^{q}T_{X}).$$
\begin{proposition}\label{deformationtorus}
Assume that there is a semiorthogonal decomposition \[ \D^{b}(X)=\langle \Ku(X),E_{1},\dots, E_{n}\rangle \] where $\{E_{1},\dots,E_{n}\}$ is an exceptional collection. 
Also assume that $\HH_{2n+1}(X)=0$ for $n\geq 1$.
Then the first-order deformation space of $J(\Ku(X))$ is
 $$H^{1}(J(\Ku(X)),T_{J(\Ku(X))})\cong \Hom(\HH_{-1}(\Ku(X)),\HH_{1}(\Ku(X))).$$
\end{proposition}

\begin{proof}
Write $V:=\HH_{1}(\Ku(X))$. Note that $V \cong \HH_{1}(X)$ by the Additivity Theorem of Hochschild Homology. Since $$\HH_{1}(\Ku(X))\cong\HH_{1}(X)\cong \sH\Omega_{1}(X)\cong\overline{\sH\Omega}_{-1}(X)\cong \overline \HH_{-1}(X)\cong \overline{\HH}_{-1}(\Ku(X)),$$ by the HKR isomorphism, there is a conjugation $\overline{V}=\HH_{-1}(\Ku(X))$. Since the tangent bundle of a torus is trivial, we have
$$H^{1}(J(\Ku(X)),T_{J(\Ku(X))})\cong H^{1}(V/\Gamma,V\otimes\mathcal{O}_{V/\Gamma}) \cong V\otimes H^{1}(V/\Gamma,\mathcal{O}_{V/\Gamma}).$$
Since $H^{1}(V/\Gamma,\mathcal{O}_{V/\Gamma})\cong\Hom_{\mathrm{anti-linear}}(V,\mathbb{C}) \cong\Hom_{\mathbb{C}}(\overline{V},\mathbb{C})$, we finally get the isomorphism
$$H^{1}(J(\Ku(X)),T_{J(\Ku(X))})\cong \Hom(\HH_{-1}(\Ku(X)),\HH_{1}(\Ku(X)))$$
as required.
\end{proof}
\par
When $\HH_{2n+1}(X)=0$, $n\geq 1$, we define a linear map from the deformations of $\Ku(X)$ to the first-order deformations of its intermediate Jacobian $J(\Ku(X))$ by the action of cohomology:
$$\HH^{2}(\Ku(X))\longrightarrow \Hom(\HH_{-1}(\Ku(X)),\HH_{1}(\Ku(X))).$$
This map can be interpreted as the derivative of the following map of ``moduli spaces''
$$\{\Ku(X)\}/\simeq  \, \longrightarrow \{J(\Ku(X))\}/\cong.$$




We now state and prove the main theorem of this section.

\begin{theorem} \label{commutative diagram theorem}
Let $X$ be a smooth projective variety. Assume $\D^{b}(X)=\langle \Ku(X), E_{1}, E_{2}, ... , E_{n}\rangle$ where $\{E_{1},E_{2}, ... , E_{n}\}$ is an exceptional collection. Then we have a commutative diagram


\[
\begin{tikzcd}[column sep=6em]
\HH^2(\Ku(X)) \arrow[r, "\gamma"]                & {\Hom(\sH\Omega_{-1}(X), \sH\Omega_1(X))} \\
\HT^2(X) \arrow[u, "\alpha'"] \arrow[ru, "\tau"] &                               \\
{H^1(X, T_X)} \arrow[u, "\mathrm{inclusion}"] \arrow[ruu]              &                              
\end{tikzcd}
\]
where $\tau$ is defined as contraction of polyvector fields.
\end{theorem}

We break the proof into several lemmas. We write $\D^{b}(X)=\langle \Ku(X), \mathcal{A}\rangle$ for the semiorthogonal decomposition, where $\mathcal{A}=\langle E_{1},E_{2}, ... ,E_{n}\rangle$. Let $P_{1}$ be the kernel of the projection to $\Ku(X)$, and $P_{2}$ the kernel of the projection to $\mathcal{A}$. 
\par
We write $\bullet$ to represent any integer. Firstly we define maps $$\alpha: \Hom_{\D^{b}(X\times X)}(\oh_\Delta,\oh_\Delta[\bullet])\rightarrow \Hom_{\D^{b}(X\times X)}(P_{1},P_{1}[\bullet])$$ and $$\beta:\Hom_{\D^{b}(X\times X)}(\oh_\Delta,\oh_\Delta\circ S_{X}[\bullet])\rightarrow \Hom_{\D^{b}(X\times X)}(P_{1},P_{1}\circ S_{X}[\bullet]).$$
There are triangles
\begin{equation} \label{kernel triangle 1}
    P_2 \longrightarrow \oh_\Delta \longrightarrow  P_1 \longrightarrow  P_2[1],
\end{equation}
\begin{equation} \label{kernel triangle 2}
    P_2 \circ S_X \longrightarrow  \oh_\Delta \circ S_X \longrightarrow  P_1 \circ S_X \longrightarrow  P_2 \circ S_X [1] .
\end{equation}
Then applying $\Hom(\mathcal{O}_{\Delta},-)$ to the triangle (\ref{kernel triangle 1}), we get a map 
\begin{equation}
\Hom_{\D^{b}(X\times X)}(\mathcal{O}_{\Delta},\mathcal{O}_{\Delta}[\bullet]) \longrightarrow \Hom_{\D^{b}(X\times X)}(\mathcal{O}_{\Delta},P_{1}[\bullet]). 
\end{equation}
 Applying $\Hom(-,P_{1})$ to the triangle (\ref{kernel triangle 1}), we get an isomorphism for any integer $\bullet$
\begin{equation}
\Hom_{\D^{b}(X\times X)}(\mathcal{O}_{\Delta},P_{1}[\bullet])\cong \Hom_{\D^{b}(X\times X)}(P_{1},P_{1}[\bullet]) ,
\end{equation}
since $\Hom_{\D^{b}(X\times X)}(P_{2},P_{1}[i])=0$ for any integer $i$ by a similar proof to \cite[ Cor 3.10]{kuznetsov2009hochschild}. 
So we get the following map $\alpha$ by composing (3) and (4):
$$\alpha: \Hom_{\D^{b}(X\times X)}(\mathcal{O}_{\Delta},\mathcal{O}_{\Delta}[\bullet]) \longrightarrow \Hom_{\D^{b}(X\times X)}(P_{1},P_{1}[\bullet]).$$
Applying $\Hom(\mathcal{O}_{\Delta},-)$ to triangle (\ref{kernel triangle 2}), we  obtain the map
\begin{equation}
\Hom_{\D^{b}(X\times X)}(\mathcal{O}_{\Delta},\mathcal{O}_{\Delta}\circ S_{X}[\bullet]) \longrightarrow \Hom_{\D^{b}(X\times X)}(\mathcal{O}_{\Delta},P_{1}\circ S_{X}[\bullet]).
\end{equation}
Similarly, applying the functor $\Hom(-,P_{1}\circ S_{X})$ to the triangle (\ref{kernel triangle 1}), we obtain an isomorphism for any integer $\bullet$
\begin{equation}
\Hom_{\D^{b}(X\times X)}(\mathcal{O}_{\Delta},P_{1}\circ S_{X}[\bullet])\cong \Hom_{\D^{b}(X\times X)}(P_{1},P_{1}\circ S_{X}[\bullet]),
\end{equation}
since $\Hom_{\D^{b}(X\times X)}(P_{2},P_{1}\circ S_{X}[i])=0$ for any integer $i$ \cite[Cor 3.10]{kuznetsov2009hochschild}. So we get the following map $\beta$ by composing (5) and (6):
$$\beta:\Hom_{\D^{b}(X\times X)}(\mathcal{O}_{\Delta},\mathcal{O}_{\Delta}\circ S_{X}[\bullet]) \longrightarrow\Hom_{\D^{b}(X\times X)}(P_{1},P_{1}\circ S_{X}[\bullet]).$$
\begin{lemma}\label{Hochschildcommutative}
There is a commutative diagram
$$\xymatrix@C=0.5cm{\Hom_{\D^{b}(X\times X)}(\mathcal{O}_{\Delta},\mathcal{O}_{\Delta}[t_{1}])\times\Hom_{\D^{b}(X\times X)}(\mathcal{O}_{\Delta},\mathcal{O}_{\Delta}\circ S_{X}[t_{2}])\ar[r]\ar[d]^{(\alpha,\beta)}&\Hom_{\D^{b}(X\times X)}(\mathcal{O}_{\Delta},\mathcal{O}_{\Delta}\circ S_{X}[t_{1}+t_{2}])\ar[d]^{\beta}\\
\Hom_{\D^{b}(X\times X)}(P_{1},P_{1}[t_{1}])
\times\Hom_{\D^{b}(X\times X)}(P_{1},P_{1}\circ S_{X}[t_{2}])\ar[r]&\Hom_{\D^{b}(X\times X)}(P_{1},P_{1}\circ S_{X}[t_{1}+t_{2}])}$$
The morphisms in the rows are the composition maps described in Definition \ref{moduleaction}.
\end{lemma}

\begin{proof}

We explain the commutative diagram. Take $t_{1}=t_{2}=0$; the general cases are similar. Let $f\in \Hom_{\D^{b}(X\times X)}(\mathcal{O}_{\Delta},\mathcal{O}_{\Delta})$ and $g\in \Hom_{\D^{b}(X\times X)}(\mathcal{O}_{\Delta},\mathcal{O}_{\Delta}\circ S_{X})$. We denote the natural morphism $\mathcal{O}_{\Delta}\rightarrow P_{1}$ by $L$. The composition $(L\otimes \identity)\circ g$ gives an element $g'$ in $$\Hom_{\D^{b}(X\times X)}(\mathcal{O}_{\Delta},P_{1}\circ S_{X})\cong \Hom_{\D^{b}(X\times X)}(P_{1},P_{1}\circ S_{X}),$$ 
that is $\beta(g)=g'$; see the commutative diagram
$$\xymatrix{\mathcal{O}_{\Delta}\ar[r]^{g}\ar[d]_{L}&\mathcal{O}_{\Delta}\circ S_{X}\ar[d]^{L\otimes \identity}\\
P_{1}\ar[r]^{g'}&P_{1}\circ S_{X}}$$
Similarly, $\alpha(f)=f'\in \Hom_{D^{b}(X\times X)}(P_{1},P_{1})$; see the commutative diagram
$$\xymatrix{\mathcal{O}_{\Delta}\ar[r]^{f}\ar[d]_{L}&\mathcal{O}_{\Delta}\ar[d]^{L\otimes \identity}\\
P_{1}\ar[r]^{f'}&P_{1}}$$
Therefore we have a commutative diagram
$$\xymatrix{\mathcal{O}_{\Delta}\ar[r]^{g}\ar[d]_{L}&\mathcal{O}_{\Delta}\circ S_{X}\ar[r]^{f\otimes \identity}\ar[d]^{L\otimes \identity}&\mathcal{O}_{\Delta}\circ S_{X}\ar[d]^{L\otimes \identity}\\
P_{1}\ar[r]^{g'}&P_{1}\circ S_{X}\ar[r]^{f'\otimes \identity}&P_{1}\circ S_{X}}$$
By the isomorphism $$\Hom_{\D^{b}(X\times X)}(\mathcal{O}_{\Delta},P_{1}\circ S_{X})\cong \Hom_{\D^{b}(X\times X)}(P_{1},P_{1}\circ S_{X}),$$ 
we have 
$$\beta((f\otimes \identity)\circ g)=(f'\otimes \identity)\circ g'.$$
\end{proof}

\begin{remark}\label{beta_is_iso}
The morphism $\beta$ is an isomorphism if $t_{1}=2$ and $t_{2}=-1$. Indeed, let $h\in \Hom_{\D^{b}(X\times X)}(\Delta_{\ast}\mathcal{O}_{X},\Delta_{\ast}\mathcal{O}_{X}\circ S_{X}[\bullet])$ for any integer $\bullet$. Define $\gamma_{P_{1}}(h)$ to be the morphism 
$$\xymatrix@C=2cm{P_{1}\circ\Delta_{\ast}\mathcal{O}_{X}\ar^{\mathrm{id} \circ h}[r]& P_{1}\circ\Delta_{\ast}\mathcal{O}_{X}\circ S_{X}[\bullet]}.$$
By \cite[Lemma 5.3]{kuznetsov2009hochschild}, there is a commutative diagram
$$\xymatrix{\Delta_{\ast}\mathcal{O}_{X}\ar[d]_{h}\ar[r]^{L}&P_{1}\ar[d]^{\gamma_{P_{1}}(h)}\\
\Delta_{\ast}\mathcal{O}_{X}\circ S_{X}\ar[r]^{L\otimes \identity}[\bullet]&P_{1}\circ S_{X}[\bullet]}$$
Hence, $\beta=\gamma_{P_{1}}$. 
Therefore, by the Theorem of Additivity \cite[Theorem 7.3]{kuznetsov2009hochschild}, $\beta$ is an isomorphism when $t_{1}=2$ and $t_{2}=-1$.
\end{remark}

Let $\alpha'$ be the following composition of maps: 
 \[ \HT^2(X) \xrightarrow{ \mathrm{IK}^{-1} } \Hom_{\D^{b}(X\times X)}(\oh_{\Delta}, \oh_{\Delta}[2]) \xrightarrow{ \ \alpha \ } \Hom_{\D^{b}(X\times X)}(P_1,P_1[2]) = \HH^2(\Ku(X)) . \]
Let $\gamma': \HT^{2}(X)\times \sH\Omega_{-1}(X)\rightarrow \sH\Omega_{1}(X)$ be the natural action of polyvectors on forms: when restricting to $H^{1}(X,T_{X})$, it is exactly the derivative of the period map.
 Define $\gamma : \HH^{2}(\Ku(X))\times\sH\Omega_{-1}(X)\rightarrow \sH\Omega_{1}(X)$ by the cohomology action as follows. Let $w\in \Hom_{\D^{b}(X\times X)}(P_{1},P_{1}[2])$. Then $\gamma(w)\colon \sH\Omega_{-1}(X)\rightarrow \sH\Omega_{1}(X)$ is defined by the commutative diagram
 $$\xymatrix{\Hom_{\D^{b}(X\times X)}(P_{1},P_{1}\circ S_{X}[1])\ar[r]^{w}&\Hom_{\D^{b}(X\times X)}(P_{1},P_{1}\circ S_{X}[1])\ar[d]^{\beta^{-1}}\\
 \Hom_{\D^{b}(X\times X)}(\mathcal{O}_{\Delta},\mathcal{O}_{\Delta}\circ S_{X}[1])\ar[u]_{\beta}&\Hom_{\D^{b}(X\times X)}(\mathcal{O}_{\Delta},\mathcal{O}_{\Delta}\circ S_{X}[1])\ar[d]^{\mathrm{IK}}\\
  \sH\Omega_{-1}(X)\ar[r]^{\gamma(w)}\ar[u]_{\mathrm{IK}^{-1}}&\sH\Omega_{-1}(X)}$$

\begin{lemma}\label{lemcommutativediag}
The following diagram is commutative:
$$\xymatrix@C6pc@R2pc{\HT^{2}(X)\times \sH\Omega_{-1}(X)\ar[r]^{\gamma'}\ar[d]^{(\alpha',\identity)}&\sH\Omega_{1}(X)\ar[d]^{\identity}\\
\HH^{2}(\Ku(X))\times \sH\Omega_{-1}(X)\ar[r]^{\gamma}&\sH\Omega_{1}(X)}$$

 \end{lemma}
 
 \begin{proof} We have a commutative diagram
 $$\xymatrix@C=0.5cm{\HT^{2}(X)\times \sH\Omega_{-1}(X)\ar[r]^{\gamma'}&\sH\Omega_{1}(X)\\
 \Hom_{\D^{b}(X\times X)}(\mathcal{O}_{\Delta},\mathcal{O}_{\Delta}[2])\times\Hom_{\D^{b}(X\times X)}(\mathcal{O}_{\Delta},\mathcal{O}_{\Delta}\circ S_{X}[-1])\ar[r]\ar[d]^{(\alpha,\beta)}\ar[u]^{\mathrm{IK}}&\Hom_{\D^{b}(X\times X)}(\mathcal{O}_{\Delta},\mathcal{O}_{\Delta}\circ S_{X}[1])\ar[d]^{\beta}\ar[u]^{\mathrm{IK}}\\
\Hom_{\D^{b}(X\times X)}(P_{1},P_{1}[2])
\times\Hom_{\D^{b}(X\times X)}(P_{1},P_{1}\circ S_{X}[-1])\ar[r]&\Hom_{\D^{b}(X\times X)}(P_{1},P_{1}\circ S_{X}[1])}$$
The upper square is commutative by \cite[Theorem 1.4]{calaque2012cualduararu}. The lower square is commutative by Lemma \ref{Hochschildcommutative}. Therefore Lemma \ref{lemcommutativediag} follows from the definitions.
 \end{proof}
  
  \begin{proof}[Proof of Theorem \ref{commutative diagram theorem}]
  By Lemma \ref{lemcommutativediag}, we obtain the commutative diagram
$$\xymatrix@C6pc@R2pc{\HH^{2}(\Ku(X))\ar[r]^{\gamma}&\Hom(\sH\Omega_{-1}(X),\sH\Omega_{1}(X))\\
\HT^{2}(X)\ar[u]^{\alpha'}\ar[ru]&\\
H^{1}(X,T_{X})\ar[u]^{\mathrm{inclusion}}\ar[ruu]_{\mathrm{d} \cP}&}$$
As a result, Theorem~\ref{commutative diagram theorem} is proved. 
\end{proof}

\begin{remark}
 In Remark \ref{beta_is_iso}, if we take $t_{2}$ to be an odd integer, then $\beta$ is also an isomorphism. The commutative diagram can be extended to the following commutative diagram:
\[
\begin{tikzcd}[column sep=6em]
\HH^2(\Ku(X)) \arrow[r, "\gamma"]                & {\bigoplus_{i}\Hom(\sH\Omega_{2i-1}(X),\sH\Omega_{2i+1}(X))} \\
\HT^2(X) \arrow[u, "\alpha'"] \arrow[ru, "\tau"] &                               \\
{H^1(X, T_X)} \arrow[u, "\mathrm{inclusion}"] \arrow[ruu]              &                              
\end{tikzcd}
\]That is to say, for any element $a\in \HH^{2}(\Ku(X))$, $\gamma(a)_i\in \Hom(\sH\Omega_{2i-1}(X),\sH\Omega_{2i+1}(X))$ is similarly defined. If $\sH\Omega_{2i+1}(X)=0$ for $i\geq 1$, then 
$\bigoplus_{i}\Hom(\sH\Omega_{2i-1}(X),\sH\Omega_{2i+1}(X))=\Hom(\sH\Omega_{-1}(X),\sH\Omega_{1}(X))$.

\end{remark}

We write $\eta$ to be the composition $\alpha'\circ (\text{inclusion})$. Note that by the Kodaira Vanishing Theorem, $h^{3,0}(X)=h^{0,3}(X)=0$ for smooth Fano threefolds. 
\begin{corollary}\label{deformationdiagram}
Let $X$ be a Fano threefold of index one or two. Note here that we have $\sH\Omega_{-1}(X)= H^{2,1}(X)$ and $\sH\Omega_{1}(X)= H^{1,2}(X)$. Then there is a commutative diagram
$$\xymatrix@C8pc@R2pc{\HH^{2}(\Ku(X))\ar[r]^{\gamma}&\Hom(H^{2,1}(X),H^{1,2}(X))\\
H^{1}(X,T_{X})\ar[ru]^{\mathrm{d} \cP}\ar[u]^{\eta}&}$$
\end{corollary}

\begin{remark}
The commutative diagram above can be regarded as the infinitesimal version of the diagram
$$\xymatrix@C4pc{\{\Ku(X)\}/\simeq\ar[r]^{\bJ}&\{J(X)\}/\cong\\
\{\mathrm{Perf}_{\mathrm{dg}}(X)\}/\cong\ar[u]\ar[ru]^{\bJ}&\\}$$
where $\mathrm{Perf}_{\mathrm{dg}}(X)$ is the dg-enhancement of perfect complexes on $X$, $\{\Ku(X)\}/\simeq$ is the Morita equivalence class of the Kuznetsov component of $X$, and $\bJ$ is the Jacobian functor
$$\bJ:\mathrm{NChow}(\mathbb{C})_{\mathbb{Q}}\rightarrow \mathrm{Ab}(\mathbb{C})_{\mathbb{Q}}$$ from the category of noncommutative Chow motives to category of abelian varieties over the complex numbers, defined in \cite[Section 1]{bernardara2013semi} and \cite[Definition 5.24]{perry2020integral}.
\end{remark}

\begin{definition} \label{infinitesimal Torelli definition}
Let $X$ be a smooth Fano threefold. Assume there is a semiorthogonal decomposition $\D^{b}(X)=\langle \Ku(X), E_{1}, E_{2}, \dots , E_{n}\rangle$ where $\{E_{1},E_{2}, \dots , E_{n}\}$ is an exceptional collection.
\begin{enumerate}
    \item The variety $X$ satisfies \emph{infinitesimal Torelli} if \[ \mathrm{d} \cP : H^1(X, T_X) \longrightarrow \Hom(\sH\Omega_{-1}(X), \sH\Omega_1(X)) \] is injective;
    \item The variety $X$ satisfies \emph{infinitesimal categorical Torelli} if the composition \[ \eta : H^1(X, T_X) \longrightarrow \HH^2(\Ku(X)) \] is injective;
    \item The Kuznetsov component $\Ku(X)$ satisfies \emph{infinitesimal Torelli} if \[ \gamma : \HH^2(\Ku(X)) \longrightarrow \Hom(\sH\Omega_{-1}(X), \sH\Omega_1(X)) \] is injective.
\end{enumerate}
\end{definition}

\begin{remark}
The definitions (2) and (3) from Definition \ref{infinitesimal Torelli definition} depend on the definition of the Kuznetsov component $\Ku(X)$ for a particular variety $X$. In other words, for a given $X$, taking the semiorthogonal complement of different exceptional collections would give different maps $\eta$.
\end{remark}

\section{Infinitesimal (categorical) Torelli theorems for Fano threefolds of index one and two} \label{infinitesimal section}
In this section, we apply Corollary~\ref{deformationdiagram} to establish infinitesimal categorical Torelli theorems for a class of prime Fano threefolds of index one and two, via the classical infinitesimal Torelli theorems for them. 

We briefly recall the classification of Fano threefolds into deformation families via their Picard rank, index, and degree (or equivalently genus). The Picard rank of a Fano threefold $X$ is the rank of $\Pic X$. In our paper, we will only be concerned with Picard rank one Fano threefolds, i.e. when $\Pic X = \mathbb{Z}$. The index of $X$ is the positive integer $i_X$ such that $K_X = -i_X H$, where $H$ is the generator of $\Pic X$. If $X$ is a Fano threefold, then $1 \leq i_X \leq 4$. The only index $4$ prime Fano threefold is $\mathbb{P}^3$ and the only index $3$ prime Fano threefold is the three dimensional quadric $Q\subset\mathbb{P}^4$. They are rigid in the sense that $H^1(X,T_X)=0$, hence infinitesimal Torelli and infinitesimal categorical Torelli theorem hold for them automatically. Thus in this paper we will be concerned with the index one and two cases. 

There are five deformation classes of Fano threefolds of index two, denoted $Y_d$. They are classified by their degree $d$, which takes the values $1 \leq d \leq 5$.

There are ten deformation classes of Fano threefolds of index one, denoted $X_{2g-2}$. They are classified by their degree (equivalently genus) $d=2g-2$, where $g$ is the genus. The allowed values of $g$ are $2 \leq g \leq 12$ and $g \neq 11$.

For a more detailed account of the Fano threefolds discussed above, see for example \cite[Section 2]{kuznetsov2003derived}.

\subsection{Classical Torelli and infinitesimal Torelli theorems}

Recall that the Torelli problem asks whether the period map $\cP:\cM\rightarrow\cD$ is injective, while the infinitesimal Torelli problem asks whether the period map has injective differential. Let $X$ be a smooth projective variety of dimension $n$ over the complex numbers $\mathbb{C}$. We say that an \emph{infinitesimal Torelli theorem} holds for $X$ if the map 
$$ \mathrm{d} \cP : H^1(X,T_X) \longrightarrow \bigoplus_{p+q=n}\mathrm{Hom}(H^p(X,\Omega_X^q),H^{p+1}(X,\Omega_X^{q-1}))$$ is injective. Besides the case of curves, the injectivity of the above map has been studied for many other varieties, which are listed in \cite[Section 1]{licht2022infinitesimal}:
\begin{itemize}
    \item hypersurfaces in (weighted) projective spaces (\cite{carlson1983infinitesimal, donagi1983generic, saito1986weak});
    \item complete intersections in projective spaces (\cite{terasoma1990infinitestimal, peters1976local, usui1976local});
    \item zero loci of sections of vector bundles (\cite{flenner1986infinitesimal});
    \item certain cyclic covers of a Hirzebruch surfaces (\cite{konno1985deformations});
    \item complete intersections in certain homogeneous K{\"a}hler manifolds (\cite{konno1986infinitesimal});
    \item weighted complete intersections (\cite{usui1978local});
    \item quasi-smooth Fano weighted hypersurfaces (\cite{fatighenti2019weighted});
    \item index one prime Fano threefolds of degree $4$ (\cite{licht2022infinitesimal}).
\end{itemize}

In particular, if $X$ is a smooth prime Fano threefold of index one or two, then the map $d\cP$ becomes
$$ \mathrm{d} \cP : H^1(X,T_X)\longrightarrow\mathrm{Hom}(H^{2,1}(X), H^{1,2}(X)),$$ since $h^{3,0}(X)=h^{0,3}(X)=0$. Then the definition of \emph{infinitesimal Torelli} for $X$ coincides with Definition~\ref{infinitesimal Torelli definition}. We have the following proposition:

\begin{proposition}[{\cite{saito1986weak, flenner1986infinitesimal}}] \label{local_Torelli_prime_Fano} \leavevmode
\begin{enumerate}
    \item Let $X_{2g-2}$ be a prime Fano threefold of index one and genus $g$. Then the map $\mathrm{d}\cP$ is injective if $g \in\{2,4,5,7\}$. In the case of $g=3$, the map $\mathrm{d}\cP$ is injective if $X$ is not hyperelliptic.
    \item Let $Y_d$ be a prime Fano threefold of index two and degree $d$. Then the map $\mathrm{d}\cP$ is injective if $d\in\{1,2,3,4,5\}$.
    \item If $X_{2g-2}$ is an index one prime Fano threefold and genus $g$ such that $g\in\{6,8,9,10,12\}$, or if $g=3$ and $X$ is hyperelliptic, then $\mathrm{d}\cP$ is not injective. 
\end{enumerate}
 
\end{proposition}

\subsection{Infinitesimal categorical Torelli theorems for Fano threefolds}

In the following sections, we study the commutative diagram constructed in Corollary \ref{deformationdiagram}, and investigate whether the infinitesimal categorical Torelli theorem defined in Definition \ref{infinitesimal Torelli definition} holds for various Fano threefolds. We use the definition of the Kuznetsov component of a Fano threefold of index one or two from the survey paper \cite{Kuznetsov_2016}, and we refer to \cite{BFT2021polyvector} and \cite{fanography} for the dimensions of 
$H^{1}(X, T_{X})$ and $H^{1}(X,\Omega^{2}_{X})$.

Let $X$ be a smooth projective variety of dimension $n$. Assume $\D^{b}(X)=\langle \Ku(X), E, \mathcal{O}_{X}\rangle$ where $E^{\vee}$ is a globally generated rank $r$ vector bundle with vanishing higher cohomology. Let $\cN_{X/\Gr}^{\vee}$ be the shifted cone lying in the triangle 
$$\cN_{X/ \Gr}^{\vee}\longrightarrow \phi^{\ast}\Omega_{\Gr(r,V)}\longrightarrow \Omega_{X},$$
where $V \coloneqq H^0(X,E^{\vee})$.

\begin{theorem}[{\cite[Theorem 8.8]{kuznetsov2009hochschild}}] \label{kuznetsov hochschild theorem 8.8}   Let the notation be as in the paragraph above. Then
\begin{enumerate}
    \item There is an exact sequence \begin{align*} \cdots &\longrightarrow \bigoplus_{p=0}^{n-1}H^{t-p}(X, \Lambda^p T_X) \longrightarrow \HH^t(\langle E, \oh_X \rangle^\bot) \longrightarrow \\ &\longrightarrow H^{t-n+2}(X, E^\bot \otimes E \otimes \omega_X^{-1}) \xrightarrow{ \ \alpha \ }  \bigoplus_{p=0}^{n-1}H^{t+1-p}(X, \Lambda^pT_X) \longrightarrow \cdots  \end{align*}
    \item There is an exact sequence \begin{align*} \cdots &\longrightarrow \bigoplus_{p=0}^{n-2}H^{t-p}(X, \Lambda^p T_X) \longrightarrow \HH^t(\langle E, \oh_X \rangle^\bot) \longrightarrow \\ &\longrightarrow H^{t-n+2}(X, \cN_{X/\Gr}^\vee \otimes \omega_X^{-1}) \xrightarrow{ \ \nu \ }  \bigoplus_{p=0}^{n-2}H^{t+1-p}(X, \Lambda^pT_X) \longrightarrow \cdots  \end{align*}
    \item If $E$ is a line bundle, then $\nu=0$ and
    $$\HH^{t}(\langle E,\mathcal{O}_{X}\rangle^\bot)\cong \bigoplus _{p=0}^{n-2}H^{t-p}(X,\Lambda^{q}T_{X})\oplus H^{t-n+2}(X,\cN_{X/\Gr}^{\vee}\otimes \omega^{-1}_{X}) . $$
\end{enumerate}
\end{theorem}

\subsection{Infinitesimal categorical Torelli for Fano threefolds of index two}

An application of part (3) of Theorem \ref{kuznetsov hochschild theorem 8.8} to the case of index two Fano threefolds of degree $d$, i.e. when $\Ku(Y_d) = \langle \oh_{Y_d}(-H), \oh_{Y_d} \rangle^\bot$, gives the following result:

\begin{theorem}[{\cite[Theorem 8.9]{kuznetsov2009hochschild}}] \label{kuznetsov hochschild theorem 8.9}
The second Hochschild cohomology of the Kuznetsov component of an index two Fano threefold of degree $d$ is given by 
\[
\HH^2(\Ku(Y_d)) = \begin{cases}
0, &   d=5 \\
k^{3}, & d=4 \\
k^{10} ,& d=3 \\
k^{20} , & d=2 \\
k^{35} , & d=1 .
\end{cases}
\]
\end{theorem}

\begin{theorem}\label{Y1234}
Let $Y_{d}$ be index two Fano threefolds of degree $1\leq d\leq 4$. 
The commutative diagrams of Corollary \ref{deformationdiagram} for $Y_{d}$ are as follows:
\begin{enumerate}
    \item $Y_{4}$: $\eta$ is an isomorphism, $\gamma$ is injective, and $d\cP$ is injective.
\[ 
\begin{tikzcd}
k^3 \arrow[r, "\gamma"]                            & k^4 \\
k^3 \arrow[u, "\eta"] \arrow[ru, "d \mathcal{P}"'] &    
\end{tikzcd}
\]

    \item $Y_{3}$: $\eta$ is an isomorphism, $\gamma$ is injective, and $d\cP$ is injective.
\[ 
\begin{tikzcd}
k^{10} \arrow[r, "\gamma"]                            & k^{25} \\
k^{10} \arrow[u, "\eta"] \arrow[ru, "d \mathcal{P}"'] &    
\end{tikzcd}
\]

   \item $Y_{2}$: $\eta$ is injective, and $d\cP$ is injective.
\[ 
\begin{tikzcd}
k^{20} \arrow[r, "\gamma"]                            & k^{100} \\
k^{19} \arrow[u, "\eta"] \arrow[ru, "d \mathcal{P}"'] &    
\end{tikzcd}
\]

   \item $Y_{1}$: $\eta$ is injective, and $d\cP$ is injective.
\[ 
\begin{tikzcd}
k^{35} \arrow[r, "\gamma"]                            & k^{441} \\
k^{34} \arrow[u, "\eta"] \arrow[ru, "d \mathcal{P}"'] &    
\end{tikzcd}
\]

\end{enumerate}
\end{theorem}

\begin{proof}


\par

The map $d\cP$ is injective for $d=1,2,3,4$ by Proposition~\ref{local_Torelli_prime_Fano}. Then by Corollary \ref{classicalandcategoricalinfinitesimal}, $\eta$ is injective. Thus for the cases $d=3,4$, the map $\eta$ is an isomorphism, hence $\gamma$ is injective.
\end{proof}

\subsection{Infinitesimal categorical Torelli for Fano threefolds of index one}

Gushel--Mukai (GM) threefolds are the index one Fano threefolds of genus $6$ (equivalently, degree $10$). They fall into two classes: ordinary and special. We study Corollary \ref{deformationdiagram} for \emph{ordinary} GM threefolds in Section \ref{ordinaryGMsubsubsection}. Finally, in Sections \ref{other_index_one_cases_subsection} and \ref{other_index_one_cases_subsection_2}, we study Corollary \ref{deformationdiagram} for the remaining index one Fano threefolds.

\begin{proposition}\label{HH2_index_one}
The second Hochschild cohomology of the Kuznetsov component of an index one Fano threefold of genus $g\geq 6$ is given by 
\[
\HH^2(\Ku(X_{2g-2})) = \begin{cases}
0, & g=12\\
k^{3}, &   g=10 \\
k^{6}, &  g=9 \\
k^{10} ,&  g=8 \\
k^{18} , & g=7 \\
k^{20} , & g=6 .
\end{cases}
\]
\end{proposition}

\begin{proof}
In the cases $g=6, 10$ and $12$, this follows from the equivalences $\Ku(X_{4d+2}) \simeq \Ku(Y_d)$ (see \cite{kuznetsov2009derived}) and Theorem \ref{kuznetsov hochschild theorem 8.9}. For the cases $g=7$ and $9$, note that we always have $\Ku(X) \simeq \D^b(C)$ for some curve $C$. The semiorthogonal decompositions are
\begin{align*}
\D^{b}(X_{16})=&\langle \D^{b}(C_{3}), \cE_{3}, \mathcal{O}_{X_{16}}\rangle\\
\D^{b}(X_{12})=&\langle \D^{b}(C_{7}), \cE_{5}, \mathcal{O}_{X_{12}}\rangle
\end{align*}
where $C_{i}$ is a curve of genus $i$ and $\cE_{r}$ is a vector bundle of rank $r$. By the $\mathrm{HKR}$ isomorphism, $\HH^{2}(C)\cong H^{2}(C,\mathcal{O}_{C})\oplus H^{1}(C,T_{C})\oplus H^{0}(C,\Lambda^{2}T_{C})= H^{1}(C,T_{C})$. The second equality follows from the fact that $C$ is of dimension one. The dimensions of $H^1(C,T_C)$ being $3i-3$ follow from a simple computations via Riemann--Roch.  For the $g=6$ case, the second Hochschild cohomology is computed in \cite[Proposition 2.12]{kuznetsov2018derived}.
\end{proof}

\subsubsection{The Gushel--Mukai case} \label{ordinaryGMsubsubsection} 
Let $X$ be an ordinary GM threefold, which is 
a quadric section of a linear section of codimension two of the Grassmannian $\mathrm{Gr}(2,5)$. There is a rank two stable vector bundle $\cE$ on $X$ which induces a morphism $\phi:X\rightarrow\mathrm{Gr}(2,5)$ such that $\cE\cong\phi^*T$, where $T$ is the tautological rank $2$ bundle on $\Gr(2,5)$.  
There is a semiorthogonal decomposition $$\D^{b}(X)=\langle \Ku(X), \cE ,\mathcal{O}_X\rangle . $$

\begin{theorem}\label{Xo10}
Let $X$ be an ordinary GM threefold. Then the commutative diagram in Corollary \ref{deformationdiagram}
is as follows:
\[ 
\begin{tikzcd}
k^{20} \arrow[r, "\gamma"]                            & k^{100} \\
k^{22} \arrow[u, "\eta"] \arrow[ru, "d \mathcal{P}"'] &    
\end{tikzcd}
\]
Moreover, $\gamma$ is injective and $\eta$ is surjective. In particular, the Kuznetsov component $\Ku(X)$ of an ordinary GM threefold satisfies infinitesimal Torelli. 
\end{theorem}

 \begin{proof}
 The moduli stack $\mathcal{X}$ of Fano threefolds of index one and degree $10$ is a $22$ dimensional smooth irreducible algebraic stack by \cite[p. 34]{debarre2012period}. 
By Section 7 of \cite{debarre2012period} we have that the differential 
\[ \mathrm{d} \cP : H^1(X, T_{X}) \longrightarrow \Hom(H^{2,1}(X), H^{1,2}(X))  \]
of the period map $ \cP : \mathcal{X} \ra \cA_{10}$ has $2$-dimensional kernel.
   Consider the kernel of $\eta$. Clearly, $\ker \eta\subset \ker d\cP$ hence $\dim \ker \eta\leq 2$.
  Since the dimension of the image of $\eta$ is less than or equal to $20$, the dimension of $\operatorname{ker}\eta$ must be $2$ and $\eta$ is surjective. Finally, since image of $d\cP$ is $20$-dimensional, so is $\gamma$, hence $\gamma$ is injective.
  \end{proof}
  
  \begin{remark}
  There is another proof for the case of ordinary GM threefolds $X$, which uses the long exact sequence in part (2) of Theorem \ref{kuznetsov hochschild theorem 8.8}:
  \begin{align*}
    0 &\longrightarrow H^0(X, \cN_{X/\Gr}^\vee(H)) \longrightarrow H^1(X, T_X) \xrightarrow{ \ \eta \ } \HH^2(\Ku(X)) \xrightarrow{ \ \nu \ } \\
    &\xrightarrow{ \ \nu \ } H^1(X, \cN_{X/\Gr}^\vee(H)) \longrightarrow H^2(X, T_X) \longrightarrow 0.
\end{align*}
  It suffices to compute $H^{0}(X,\cN^{\vee}_{X/\Gr}(1))$ and $H^{1}(X,\cN^{\vee}_{X/\Gr}(1))$. Since $\cN_{X/\Gr}$ is the restriction of $\mathcal{O}_{\Gr}(1)\oplus \mathcal{O}_{\Gr}(1)\oplus \mathcal{O}_{\Gr}(2)$ to $X$, we have that $\cN_{X/\Gr}^{\vee}(1)$ is the restriction of $\mathcal{O}_{\Gr}\oplus \mathcal{O}_{\Gr}\oplus \mathcal{O}_{\Gr}(-1)$ to $X$, which is $\mathcal{O}_{X}\oplus \mathcal{O}_{X}\oplus \mathcal{O}_{X}(-H)$. Then $H^{0}(X,\cN^{\vee}_{X/\Gr}(H))=k^{2}$ and $H^{1}(X,\cN^{\vee}_{X/\Gr}(H))=0$ by the Kodaira Vanishing Theorem. Thus, $\eta$ is surjective with two dimensional kernel, hence $\gamma$ is injective.
  \end{remark}




\subsubsection{The cases of $X_{18}$, $X_{16}$, $X_{14}$, and $X_{12}$} \label{other_index_one_cases_subsection}
Consider the index one prime Fano threefolds of genus $7\leq g\leq 10$. They are \begin{enumerate}
    \item $X_{12}$, $g=7$: a linear section of a connected component of the orthogonal Lagrangian Grassmannian $\mathrm{OGr}_+(5,10)\subset \mathbb{P}^{15}$;
    
    \item $X_{14}$, $g=8$: a linear section of $\mathrm{Gr}(2,6)\subset \mathbb{P}^{14}$;
    
    \item $X_{16}$, $g=9$: a linear section of the Lagrangian Grassmannian $\mathrm{LGr}(3,6)\subset \mathbb{P}^{13}$;
    
    \item $X_{18}$, $g=10$: a linear section of the homogeneous space $G_2/P\subset \mathbb{P}^{13}$.
\end{enumerate}

\begin{theorem}\label{derivedinvariantcommutative}
Let $X$ and $Y$ be smooth prime Fano threefolds of index one or index two. Suppose there is an equivalence of their Kuznetsov components $\Ku(X)\simeq \Ku(Y)$ which is induced by a Fourier--Mukai functor. Then there is a commutative diagram 
$$\begin{tikzcd}
  \HH^{2}(\Ku(X))\arrow[r, "\gamma_{X}"]\arrow[d, "\cong"]&\Hom(\HH_{-1}(\Ku(X))), \HH_{1}(\Ku(X)))\arrow[d, "\cong"]\\
  \HH^{2}(\Ku(Y))\arrow[r, "\gamma_{Y}"]&\Hom(\HH_{-1}(\Ku(Y)), \HH_{1}(\Ku(Y)))
  \end{tikzcd}$$
\end{theorem}

\begin{proof}
By Theorem \ref{derivedinvariant}, there is a commutative diagram
$$\begin{tikzcd}
\HH^{2}(\Ku(X))\times \HH_{-1}(\Ku(X))\arrow[r,"\gamma_{X}"]\arrow[d,"\cong"]& \HH_{1}(\Ku(X))\arrow[d,"\cong"]\\
\HH^{2}(\Ku(Y))\times \HH_{-1}(\Ku(Y))\arrow[r, "\gamma_{Y}"]& \HH_{1}(\Ku(Y))
\end{tikzcd}$$
The maps in the rows are defined as the cohomology action on homology. Hence the commutative diagram in the theorem follows.
\end{proof}

\begin{remark}
When $X$ and $Y$ are Fano threefolds of index one and two, respectively, we have $\HH_{-1}(X)\cong H^{2,1}(X)$ and $\HH_{-1}(Y)\cong H^{2,1}(Y)$, respectively. Hence we obtain a commutative diagram
 $$\begin{tikzcd}
  \HH^{2}(\Ku(X))\arrow[r, "\gamma_{X}"]\arrow[d, "\cong"]&\Hom(H^{2,1}(X), H^{1,2}(X))\arrow[d, "\cong"]\\
  \HH^{2}(\Ku(Y))\arrow[r, "\gamma_{Y}"]&\Hom(H^{2,1}(Y), H^{1,2}(Y))
  \end{tikzcd}$$
 where $\gamma_{X}$ and $\gamma_{Y}$ are the maps constructed in Theorem \ref{deformationdiagram}.
\end{remark}

\begin{theorem} \label{X12141618}
The diagrams in Corollary \ref{deformationdiagram} for $X_{18}$, $X_{16}$, $X_{14}$, and $X_{12}$ are as follows:
\begin{enumerate}
    \item $X_{18}$: $\gamma$ is injective.
\[ 
\begin{tikzcd}
k^{3} \arrow[r, "\gamma"]                            & k^{4} \\
k^{10} \arrow[u, "\eta"] \arrow[ru, "d \mathcal{P}"'] &    
\end{tikzcd}
\]
   \item $X_{16}$: $\gamma$ is injective.
\[ 
\begin{tikzcd}
k^{6} \arrow[r, "\gamma"]                            & k^{9} \\
k^{12} \arrow[u, "\eta"] \arrow[ru, "d \mathcal{P}"'] &    
\end{tikzcd}
\]
    \item $X_{14}$: $\gamma$ is injective.
\[ 
\begin{tikzcd}
k^{10} \arrow[r, "\gamma"]                            & k^{25} \\
k^{15} \arrow[u, "\eta"] \arrow[ru, "d \mathcal{P}"'] &    
\end{tikzcd}
\]
   \item $X_{12}$: $\gamma$ is injective, $d\cP$ is injective, and $\eta$ is an isomorphism.
\[ 
\begin{tikzcd}
k^{18} \arrow[r, "\gamma"]                            & k^{49} \\
k^{18} \arrow[u, "\eta"] \arrow[ru, "d \mathcal{P}"'] &    
\end{tikzcd}
\]
\end{enumerate}
\end{theorem}

\begin{proof}
In the cases of $X_{18}$, $X_{16}$, and $X_{12}$ we always have $\Ku(X)\simeq \D^{b}(C)$ for some curve $C$. The semiorthogonal decompositions are
\begin{align*}
\D^{b}(X_{18})=&\langle \D^{b}(C_{2}), \cE_{2}, \mathcal{O}_{X_{18}}\rangle\\
\D^{b}(X_{16})=&\langle \D^{b}(C_{3}), \cE_{3}, \mathcal{O}_{X_{16}}\rangle\\
\D^{b}(X_{12})=&\langle \D^{b}(C_{7}), \cE_{5}, \mathcal{O}_{X_{12}}\rangle
\end{align*}
  where $C_{i}$ is a curve of genus $i$ and $\cE_{r}$ is a vector bundle of rank $r$.
 We write $X$ for $X_{18}$, $X_{16}$, and $X_{12}$. By the $\mathrm{HKR}$ isomorphism, $\HH^{2}(C)\cong H^{2}(C,\mathcal{O}_{C})\oplus H^{1}(C,T_{C})\oplus H^{0}(C,\wedge^{2}T_{C})= H^{1}(C,T_{C})$.
  Note that we always refer to the version of $\mathrm{HKR}$ twisted by $\mathrm{IK}$, as the ``twisted $\mathrm{HKR}$". This $\mathrm{IK}$ isomorphism preserves the module structure, where the geometric side is the action of polyvector fields on differential forms. Thus by Theorem \ref{derivedinvariantcommutative} there is a commutative diagram
  $$
  \begin{tikzcd}
  \HH^{2}(\Ku(X))\arrow[r, "\gamma"]\arrow[d, "\cong"]&\Hom(H^{2,1}(X), H^{1,2}(X))\arrow[d, "\cong"]\\
  H^{1}(C,T_{C})\arrow[r, "d\cP_{C}"]& \Hom(H^{1,0}(C), H^{0,1}(C))
  \end{tikzcd}$$
  Therefore $\gamma$ is injective for each $X$ since $d\cP_{C}$ is injective for each $C:=C_i$. Indeed, $C_3$ is a plane quartic curve (\cite[Section 3.1]{brambilla2009rank}), which is a canonical curve in $\mathbb{P}^2$. Similarly, $C_7$ is also a canonical curve in $\mathbb{P}^6$ by \cite[Section 1]{iliev2007parametrization}. Thus they are both non-hyperelliptic. 
  
  \par
  
  For the case $X_{14}$, it is known that $\Ku(X_{14})\simeq \Ku(Y_{3})$ by \cite{kuznetsov2009derived} (so also their Hochschild cohomologies are the same). Then by Theorem \ref{derivedinvariantcommutative} there is a commutative diagram
  $$\begin{tikzcd}
  \HH^{2}(\Ku(X_{14}))\arrow[r, "\gamma_{X_{14}}"]\arrow[d, "\cong"]&\Hom(H^{2,1}(X_{14}), H^{1,2}(X_{14}))\arrow[d, "\cong"]\\
  \HH^{2}(\Ku(Y_{3}))\arrow[r, "\gamma_{Y_{3}}"]&\Hom(H^{2,1}(Y_{3}), H^{1,2}(Y_{3}))
  \end{tikzcd}$$
   Then $\gamma_{X_{14}}$ is injective since $\gamma_{Y_{3}}$ is injective, by Theorem \ref{Y1234}.
\end{proof}

\subsubsection{The cases of $X_{8}$, $X_6$, non hyperelliptic $X_{4}$, and $X_{2}$} \label{other_index_one_cases_subsection_2}
Consider the index one prime Fano threefolds of genus $2\leq g\leq 5$. They are
\begin{enumerate}
     \item $X_2$, $g=2$: a double cover of $\mathbb{P}^3$ branched in a surface of degree six;
    
    \item $X_4$, $g=3$: either a quartic threefold, or the double cover of a smooth quadric threefold branched in an intersection with a quartic;
    
    \item $X_6$, $g=4$: a complete intersection of a cubic and a quadric;
    
    \item $X_8$, $g=5$: a complete intersection of three quadrics. 
    
\end{enumerate}

In these cases, the Kuznetsov components are defined as $\langle \mathcal{O}_X\rangle^{\perp}$ by \cite[Definition 6.5]{bayer2017stability}. 

\begin{theorem} \label{X248}
The diagrams in Corollary \ref{deformationdiagram} for $X_{8}$, $X_6$, $X_{4}$, and $X_{2}$ are as follows:

\begin{enumerate}
    \item $X_{8}$: $\gamma$ is injective, $\eta$ is an isomorphism, and $d\cP$ is injective.
\[ 
\begin{tikzcd}
k^{27} \arrow[r, "\gamma"]                            & k^{196} \\
k^{27} \arrow[u, "\eta"] \arrow[ru, "d \mathcal{P}"'] &    
\end{tikzcd}
\]

 \item $X_{6}$: $\gamma$ is injective, $\eta$ is an isomorphism, and $d\cP$ is injective.
\[ 
\begin{tikzcd}
k^{34} \arrow[r, "\gamma"]                            & k^{400} \\
k^{34} \arrow[u, "\eta"] \arrow[ru, "d \mathcal{P}"'] &    
\end{tikzcd}
\]

    \item \begin{enumerate}
        \item If $X_4$ is a smooth quartic threefold, then $\gamma$ is injective, $\eta$ is an isomorphism, and $d\cP$ is injective.
        \item If $X_4$ is a hyperelliptic Fano threefold, then $\eta$ is an isomorphism and neither of $\gamma$ and $d\cP$ is injective.
    \end{enumerate}
\[ 
\begin{tikzcd}
k^{45} \arrow[r, "\gamma"]                            & k^{900} \\
k^{45} \arrow[u, "\eta"] \arrow[ru, "d \mathcal{P}"'] &    
\end{tikzcd}
\]
    \item $X_{2}$: $\gamma$ is injective, $\eta$ is an isomorphism, and $d\cP$ is injective.
\[ 
\begin{tikzcd}
k^{68} \arrow[r, "\gamma"]                            & k^{2704} \\
k^{68} \arrow[u, "\eta"] \arrow[ru, "d \mathcal{P}"'] &    
\end{tikzcd}
\]
\end{enumerate}
\end{theorem}

\begin{proof}
First, we prove that $\eta$ is an isomorphism in each case. We write $X$ for $X_{8}$, $X_6$, $X_{4}$, and $X_{2}$. Note that $\D^b(X)=\langle \Ku(X), \mathcal{O}_{X}\rangle$. Denote by $P_{1}$ the kernel of the left projection to $\Ku(X)$, and $P_{2}$ the kernel of the right projection to $\langle \mathcal{O}_{X}\rangle$.
There is a triangle
$$P_{2}\rightarrow \Delta_{\ast}\mathcal{O}_{X}\rightarrow P_{1}\rightarrow P_{2}[1] . $$
Applying the functor $\Delta^{!}$ to the triangle, we obtain the diagram
$$\begin{tikzcd}
\Delta^{!}P_{2}\arrow[r]\arrow[d,"\cong"]&\Delta^{!}\Delta_{\ast}\mathcal{O}_{X}\arrow[r,"L"]\arrow[d,"\cong"]& \Delta^{!}P_{1}\arrow[d,"\mathrm{id}"]\\
\omega^{-1}_{X}[-3]\arrow[r,"w"]&\bigoplus^{3}_{p=0}\Lambda^{p}T_{X}[-p]\arrow[r]&\Delta^{!}P_{1} 
\end{tikzcd}$$
By \cite[Theorem 8.5]{kuznetsov2009hochschild}, the map $w$ is an isomorphism onto the third summand.
Applying $\Hom^{2}(\mathcal{O}_{X},-)$, we obtain the commutative diagram
$$\begin{tikzcd}
\Hom_{\D^{b}(X\times X)}(\mathcal{O}_{X},\Delta^{!}\Delta_{\ast}\mathcal{O}_{X}[2])\arrow[d,"\cong"]\arrow[r,"L"]& \Hom_{\D^{b}(X\times X)}(\mathcal{O}_{X}, \Delta^{!}P_{1}[2])\arrow[d,"\mathrm{id}"]\\
\Hom_{\D^{b}(X\times X)}(\mathcal{O}_{X},\bigoplus^{3}_{p=0}\Lambda^{p}T_{X}[2-p])\arrow[r]&\Hom_{\D^{b}(X\times X)}(\mathcal{O}_{X}, \Delta^{!}P_{1}[2])\\
\Hom_{\D^{b}(X\times X)}(\mathcal{O}_{X}, \bigoplus^{2}_{p=0}\Lambda^{p}T_{X}[2-p])\arrow[ru,"\cong"]\arrow[u,"\cong"]& 
\end{tikzcd}$$
Thus, the morphism $L$ is an isomorphism. However, $L$ is naturally isomorphic to the morphism
$$\Hom_{\D^{b}(X\times X)}(\Delta_{\ast}\mathcal{O}_{X}, \Delta_{\ast}\mathcal{O}_{X}[2])\rightarrow \Hom_{\D^{b}(X\times X)}(\Delta_{\ast}\mathcal{O}_{X}, P_{1}[2])\cong \Hom_{\D^{b}(X\times X)}(P_{1}, P_{2}[2]).$$
That is to say the map $\alpha': \HT^{2}(X)\rightarrow \HH^{2}(\Ku(X))$ constructed in Theorem \ref{commutative diagram theorem} is an isomorphism. By \cite[Appendix A]{BFT2021polyvector}, $H^{0}(X, \Lambda^{2}T_{X})=0$ and it is clear that  $H^{2}(X,\mathcal{O}_{X})=0$. Then the inclusion $H^1(X,T_X)\hookrightarrow \HT^2(X)=H^{0}(X, \Lambda^{2}T_{X})\oplus H^1(X,T_X)\oplus H^2(X,\oh_X)=H^1(X,T_X)$ is an isomorphism. Thus $\eta$ (given as the composition of the inclusion and $\alpha'$) is an isomorphism. 

The map $d\cP$ is injective for $X_{8}$, $X_6$, non-hyperelliptic $X_{4}$, and $X_2$ by Proposition~\ref{local_Torelli_prime_Fano}. 
Then $\gamma$ is injective for these cases because $\eta$ is an isomorphism. The map $d\cP$ is not injective for hyperelliptic $X_4$, thus $\gamma$ is not injective. 
\end{proof}

\begin{remark}
\label{alternative_version_Ku_6}
In the literature \cite[Example 3.17]{kuznetsov2021semiorthogonal}, one can define an alternative version of the Kuznetsov component of $X_6$ as $\Ku(X):=\langle \cU_1,\oh_X\rangle^{\perp}$, where $\cU_1$ is restriction of one of the spinor bundles on the quadric $M\subset\mathbb{P}^5$. Then we still have the commutative diagram from Corollary~\ref{deformationdiagram} as follows:
\[ 
\begin{tikzcd}
\mathrm{HH}^2(\Ku(X)) \arrow[r, "\gamma"]                            & k^{400} \\
k^{34} \arrow[u, "\eta"] \arrow[ru, "d \mathcal{P}"'] &    
\end{tikzcd}
\]
Since $d\cP$ is injective, $\eta$ is injective. Thus infinitesimal categorical Torelli holds for $X_6$. By Lemma~\ref{2nd_HH_6}, the second Hochschild cohomology is given by $\mathrm{HH}^2(\Ku(X))\cong k^{34}$, thus $\eta$ is an isomorphism, hence $\gamma$ is injective.


\begin{lemma}
\label{2nd_HH_6}
Let $X$ be an index one prime Fano threefold of genus $g=4$ (degree $6$) with the semiorthogonal decomposition 
$$\D^b(X)=\langle\Ku(X),\cU_1,\oh_X\rangle.$$
Then $\mathrm{HH}^2(\Ku(X))\cong H^1(X,T_X)\cong k^{34}$. 
\end{lemma}

\begin{proof}
Consider the tautological short exact sequence
$$0\longrightarrow \cU_1\longrightarrow\oh_X^{\oplus 4}\longrightarrow\cQ_1\longrightarrow 0,$$
where $\cQ_1$ the tautological quotient bundle. It is known that $\cQ_1\cong\cU^{\vee}_2$ by \cite[Theorem 2.8(ii)]{ottaviani1988spinor}. It is easy to see that $\cU^{\perp}_1\cong\cQ_1^{\vee}\cong\cU_2$. By Theorem~\ref{kuznetsov hochschild theorem 8.8}, there is a long exact sequence
$$\dots\longrightarrow H^0(X,\cU_1^{\perp}\otimes\cU_1^{\vee})\longrightarrow H^1(X,T_X)\longrightarrow\mathrm{HH}^2(\Ku(X))\longrightarrow H^1(X,\cU_1^{\perp}\otimes\cU_1^{\vee})\longrightarrow\cdots.$$
Note that $H^{\bullet}(X,\cU_1^{\perp}\otimes\cU_1^{\vee})\cong\mathrm{Hom}^{\bullet}(\cU_1,\cU_2)=0$ since $\cU_1,\cU_2$ are completely orthogonal by \cite[Remark 5.19]{kuznetsov2021semiorthogonal}. Hence $\mathrm{HH}^2(\Ku(X))\cong H^1(X,T_X)\cong k^{34}.$
\end{proof}

\end{remark}

\bibliographystyle{alpha}
{\small{\bibliography{hochschild}}}

\end{document}